\documentclass[final, times, 14p]{elsarticle}
\usepackage{amsmath}
\usepackage{amssymb}
\usepackage{bm}
\usepackage{graphicx}
\usepackage{longtable} 
\usepackage{array}     
\usepackage{lscape}    
\usepackage{multirow}
\usepackage{xcolor}
\usepackage{subcaption}
\usepackage{algorithm}
\usepackage{algpseudocode}

\newcommand{\strainrate}{\bm{\dot\varepsilon}}

\usepackage{booktabs}           
\usepackage{siunitx}            
\usepackage{amsmath}            
\usepackage{multirow}           

\usepackage{soul}
\usepackage{xcolor}
\definecolor{darkblue}{RGB}{0,0,139}  
\definecolor{darkred}{RGB}{139,0,0}  

\usepackage{makecell}

\newcommand{\q}[1]{``#1''}
\begin{document}
	
\title{Variational Derivation of the Weak Form for Viscous Flows
	Involving Free Surfaces with Implications for Geodynamic Simulations}

\author[1]{Karim Norouzi Moghanjoghi\corref{cor1}}
\ead{norouzim@uni-bremen.de}

\author[1]{Javier García-Pintado}
\ead{jgarciapintado@marum.de}

\author[2]{Daniel R. Shapero}
\ead{shapero@uw.edu}

\author[1]{Marta Perez-Gussinye}
\ead{mpgussinye@marum.de}

\cortext[cor1]{Corresponding author}

\affiliation[1]{
	organization={Department of Geosciences, MARUM -- Center for Marine
		Environmental Sciences, University of Bremen},
	city={Bremen},
	country={Germany}
}

\affiliation[2]{
	organization={Polar Science Center, Applied Physics Laboratory,
		University of Washington},
	city={Seattle, WA},
	country={USA}
}
	\begin{abstract}
		Geodynamics simulations involve solving the Stokes equations, where boundary conditions play a crucial role in governing the internal dynamics and surface topography evolution. A common approach is to model the interface in contact with air as a free surface. This interface is therefore stress-free and deforms in response to internal force balances. This has proven essential for understanding plate tectonics, from the evolution of rifted margins to subduction processes. However, using the free surface directly can introduce instabilities and oscillations in the numerical solution, commonly referred to as the 'drunken sailor' instability, which forces the model to use small time steps to maintain stability. To mitigate this issue, geodynamics models often include a stabilization term that takes into account the effect of gravity forces due to free surface movement. So far, the stabilization term for the free surface has been derived primarily on a heuristic basis, offering a practical solution while leaving room for a more rigorous mathematical foundation.
		In this study, we apply variational calculus to the energy minimization formulation of incompressible Stokes flow. By using the Gateaux derivative, we derive the complete weak form of the Stokes equations and demonstrate that stabilization terms are a natural consequence of the free surface assumption. We only consider the linear form of viscosity, as has been used in the literature. We find that there are two terms that appear as a boundary integral in the weak form. One term is due to the internal viscous deformation and the other is due to gravitational body force. Two test cases are solved with the new weak formulation to illustrate the consistency of the previously suggested stabilization terms with the newly proposed derivation. Our results mathematically justify and complement the application of the heuristically derived gravity stabilization.
	\end{abstract}
	\maketitle
	
	\begin{highlights}
		\item Free surface assumption is systematically explored for linear Stokes equation.
		\item The complete weak form is derived for problems with free surface.
		\item Weak form is derived from minimization of the energy integral of Stokes flow.
		\item The free surface and its relation to Neumann boundary is revealed.
		\item The approach provides a mathematical basis for utilization of gravity stabilization term.
	\end{highlights}
	
	\section{Introduction}
	
	Geodynamic simulations employ viscous flow for modeling the deformation in the crust and mantle. The type of boundary conditions for these simulations is of great significance, with effects ranging from the development of an adequate topography to sub-lithospheric dynamics. In geodynamics problems, the top surface of the domain is in contact with air, and therefore that boundary is not constrained by external forces. This boundary, called the free surface, moves as a consequence of the balance between internal stresses and body forces, which determines its shape and motion. For example, it has been shown that both the timing and behavior of subducting slabs are controlled by the free surface boundary
	condition. The difference between a free surface model and a free-slip boundary is considerable.  \cite{kausStabilizationAlgorithmGeodynamic2010, schmelingBenchmarkComparisonSpontaneous2008}.  In the comparative study by \citet{schmelingBenchmarkComparisonSpontaneous2008}, where most of the participating codes used the "sticky air" assumption, it was concluded that the parametrization of the sticky air layer controls the shape and timing of slab subduction. The importance of boundary condition choices, particularly the free surface, was also highlighted by   \citet{quinquisRoleBoundaryConditions2011}, who demonstrated that slab subduction is only possible if a true free surface is used or when a layer of low viscosity and low-density material exists between the lithosphere and a free-slip boundary.
	
	However, including the free surface in geodynamic models limits the time step size. With larger time steps, the velocity field starts to oscillate around the true solution, which eventually leads to numerical blow-up, known as the "drunken sailor" effect \cite{kausStabilizationAlgorithmGeodynamic2010, kramerImplicitFreeSurface2012}. To avoid this effect, some researchers have used what is called the "sticky air" approximation,  which assumes a viscous layer above the top surface \citep{crameriComparisonNumericalSurface2012}. This approach has been shown to be sufficiently accurate for certain combinations of viscosity, sticky air depth, and characteristic length scale. However, the Stokes equation becomes much more difficult to solve as the difference between viscosity of air and rock is very large, therefore, a non-realistic viscosity is assumed for air which is many orders of magnitude greater than actual viscosity of air. In finite difference grids, \citet{duretzFreeSurfaceCapturing2016} performed interface capturing with slight modifications to the computational stencil, demonstrating first-degree accuracy. Implicit methods have also been explored to mitigate the drunken sailor effect \citep{furuichiImplicitSolutionMaterial2015, popovSLIM3DToolThreedimensional2008, kramerImplicitFreeSurface2012}.  
	
	Other efforts have been directed to model a true free surface by applying stabilization forces to compensate for the "drunken sailor" effect, thereby allowing larger time steps to be used. Various stabilization methods have been proposed on a heuristic basis \cite{kausStabilizationAlgorithmGeodynamic2010, kramerImplicitFreeSurface2012, andres-martinezNewFreesurfaceStabilization2015, roseStabilityAccuracyFree2017}. For example, \citet{kausStabilizationAlgorithmGeodynamic2010} uses a Taylor series expansion of the velocity field components, resulting in a stabilization term involving gravitational forces and a parameter $\theta$ to scale the time step and control method accuracy. Another implicit method  \citet{kramerImplicitFreeSurface2012} introduces the free surface elevation as a new variable in the system of equations. Using the same kinematic boundary conditions on the free surface \citet{ahlkronaFullyImplicitSecond2026} developed a fully implicit method for free surface simulation. The later methods solve an additional equation besides the momentum and mass conservation. A penalty approach for the free surface, incorporating various slope configurations and coupled with a penalty Stokes solver, was suggested by \citet{andres-martinezNewFreesurfaceStabilization2015}. A non-standard finite difference method proposed by \citet{roseStabilityAccuracyFree2017} showed better stability compared to the quasi-implicit method of \citet{kausStabilizationAlgorithmGeodynamic2010}. An extensive study on the accuracy of both non-standard and quasi-implicit methods was conducted in \citep{roseStabilityAccuracyFree2017}. These methods, though effective and easy to use, are not based on a first principles mathematical derivation. \\

	In this study, we consider the Stokes flow problem in a domain that moves and deforms over time due to the existence of a free surface. We use the internal energy minimization formulation of Stokes flow proposed by \cite{boffiMixedFiniteElement2013}, which allows us to derive the complete form of the forces acting on the free surface. We present a new weak form that includes a viscous and gravity stabilization terms, where the latter coincides with that proposed by \cite{kausStabilizationAlgorithmGeodynamic2010}. Thus, we complement the accepted and widely used stabilization term as part of a more general weak form arising from spatial variation of a deforming domain.
	
	\section{Governing equations}

	For the sake of completeness, here and in the next section, we reproduce the governing differential equations and the standard weak form of the Stokes equations used in geodynamics simulations.
	Let $\Omega \subset \mathbb{R}^d$ be a bounded domain with boundary $\partial\Omega$. The unknowns are velocity field $\mathbf{u} : \Omega \to \mathbb{R}^d$ and the pressure field $p : \Omega \to \mathbb{R}$. The incompressible Stokes equations are
	\begin{align}
		\nabla \cdot \boldsymbol{\sigma} + \mathbf{f} &= \mathbf{0} \quad && \text{in }  \Omega \\
		\nabla \cdot \mathbf{u} &= 0 \quad  && \text{in }  \Omega
	\end{align}
	where $\mathbf{f}$ is the body force (usually gravity), and $\boldsymbol{\sigma}$ is the Cauchy stress tensor. The stress tensor is decomposed into its deviatoric and volumetric (pressure) parts
	\begin{equation}
		\boldsymbol{\sigma} = \bm{\sigma_d} -p \mathbf{I} = 2\mu \, \strainrate(\mathbf{u}) -p \mathbf{I},
		\label{stress_tensor_def}
	\end{equation}
	where $\bm{\sigma_d}$ is the deviatoric stress tensor and
	\begin{equation}
		\strainrate(\mathbf{u}) = \frac{1}{2} \left( \nabla \mathbf{u} + (\nabla \mathbf{u})^\top \right)
	\end{equation}
	is the symmetric gradient (strain rate tensor), and $\mu > 0$ the dynamic effective viscosity.
	
	Thus, in terms of unknown variables, the strong form of the incompressible Stokes flow on the domain
	$\Omega \subset \mathbb{R}^d$ with boundary $\Gamma = \Gamma_\mathbf{D} \cup \Gamma_\mathbf{N}$ is given by
	
	\begin{align}
		\nabla \cdot 2 \mu \, \strainrate(\mathbf{u})
		- \nabla p
		+ \mathbf{f}
		&= \mathbf{0}
		&& \text{in } \Omega,
		\label{eq:momentum_balance}
		\\
		\nabla \cdot \mathbf{u}
		&= 0
		&& \text{in } \Omega,
		\label{eq:divergence_free}
		\\
		\mathbf{u}
		&= \mathbf{u}_0
		&& \text{on } \Gamma_D,
		\label{eq:dirichlet_bc}
		\\
		\boldsymbol{\sigma}  \mathbf{n}
		&= \mathbf{t}
		&& \text{on } \Gamma_N,
		\label{eq:neumann_bc}
	\end{align}
	
	where $\Gamma_N$ is non-empty. Equation \eqref{eq:momentum_balance} is the conservation of the linear momentum and equation \eqref{eq:divergence_free} is the conservation of mass.
	$\Gamma_\mathbf{D}$ represents a Dirichlet boundary where velocity is specified, and $\Gamma_\mathbf{N}$ represents a Neumann type boundary, where a traction vector $\mathbf{t}$ is specified. For the sake of  brevity, we assume the traction force is $\mathbf{t}$ zero in this study. It should be noted that the problem contains a Neumann boundary $\Gamma_\mathbf{N}$ and therefore there is no pressure null space present. 
	
	\section{Derivation of the standard weak form}
	
	Here we provide the derivation of the weak form of the standard Stokes flow, based on the multiplication of the PDE terms with test functions and integrating over the domain. Using finite elements to find the solution to a partial differential equation starts with forming a \emph{weighted residual} (or \emph{weak form}).
	We define the Lebesgue and Sobolev spaces used in this work as follows.
	The space of square-integrable scalar functions is given by
	\begin{equation}
		L^2(\Omega)
		:=
		\left\{
		p : \Omega \to \mathbb{R}
		\;\middle|\;
		\int_\Omega |p(x)|^2 \, dx < \infty
		\right\},
		\label{eq:L2}
	\end{equation}
	while the Sobolev space of vector-valued functions $\mathbf u$ with square-integrable weak
	gradients is defined as
	\begin{equation}
		H^1(\Omega)^d
		:=
		\left\{
		\mathbf{u} \in L^2(\Omega)^d
		\;\middle|\;
		\nabla \mathbf{u} \in L^2(\Omega)^{d \times d}
		\right\}.
		\label{eq:H1}
	\end{equation}
	where $d$ is the spatial dimension of the domain $\Omega$.
	Let $\mathbf{V_0}$ be the space of admissible velocity test functions defined as:
	
	\begin{equation}
		\mathbf{V}_0
		=
		H^1_{0,\Gamma_D}(\Omega)^d
		=
		\left\{
		\mathbf{v} \in H^1(\Omega)^d
		\;\middle|\;
		\mathbf{v} = \mathbf{0} \text{ on } \Gamma_D
		\right\}.
		\label{eq:velocity_test_space}
	\end{equation}
	and the corresponding trial function space defined as:
	\begin{equation}
		\mathbf{V}
		=
		\left\{
		\mathbf{u} \in H^1(\Omega)^d
		\;\middle|\;
		\mathbf{u} = \mathbf{u}_D \text{ on } \Gamma_D
		\right\}.
		\label{eq:velocity_trial_space}
	\end{equation}
	
	Let $\mathbf{V}_0$ and $Q = L^2(\Omega)$ be the spaces for velocity and pressure test functions. Let $\mathbf{v} \in \mathbf{V}_0$ and $q \in Q$ be arbitrary test functions.
	Multiplying the momentum equation \eqref{eq:momentum_balance} by $\mathbf{v}$ and integrating over $\Omega$ gives
	
	\begin{equation}
		\int_\Omega \mathbf{v} \cdot \left(\nabla \cdot \boldsymbol{\sigma} \right) \, d\mathbf{x} + \int_\Omega \mathbf{v} \cdot \mathbf{f} \, d\mathbf{x} = 0.
	\end{equation}
	
	By applying integration by parts to the divergence term in the momentum equation
	\[
	\int_\Omega \mathbf{v} \cdot (\nabla \cdot \boldsymbol{\sigma}) \, d\mathbf{x} = -\int_\Omega \boldsymbol{\sigma} : \nabla \mathbf{v} \, d\mathbf{x} + \int_\Omega  \nabla \cdot (\mathbf{v} \cdot \boldsymbol{\sigma}) \, d\mathbf{x},
	\]
	and the divergence theorem
	\[
	\int_\Omega  \nabla \cdot (\mathbf{v} \cdot \boldsymbol{\sigma}) \, d\mathbf{x} =  \int_{\Gamma_N} (\boldsymbol{\sigma} \cdot \mathbf{n}) \cdot \mathbf{v} \, dS,
	\]
	we obtain:
	
	\begin{equation}
		\int_\Omega \boldsymbol{\sigma} : \nabla \mathbf{v} \, d\mathbf{x} - \int_\Omega \mathbf{f} \cdot \mathbf{v} \, d\mathbf{x} -  \int_{\Gamma_N} (\boldsymbol{\sigma} \cdot \mathbf{n}) \cdot \mathbf{v} \, dS = 0.
	\end{equation}
	Here the colon denotes the double contraction, Frobenius product, of tensors.
	Substituting with the definition of the stress tensor given in \eqref{stress_tensor_def}:
	
	\begin{align}
		\int_\Omega (2\mu \, \strainrate(\mathbf{u}) -p\mathbf{I}) : \nabla \mathbf{v} \, d\mathbf{x} - \int_\Omega \mathbf{f} \cdot \mathbf{v} \, d\mathbf{x} - \int_{\Gamma_N} (\boldsymbol{\sigma} \cdot \mathbf{n}) \cdot \mathbf{v} \, dS  &= 0 \\
		\Rightarrow \int_\Omega 2\mu \, \strainrate(\mathbf{u}) : \nabla \mathbf{v} \, d\mathbf{x} - \int_\Omega p \, \nabla \cdot \mathbf{v} \, d\mathbf{x} - \int_\Omega \mathbf{f} \cdot \mathbf{v} \, d\mathbf{x} - \int_{\Gamma_N} (\boldsymbol{\sigma} \cdot \mathbf{n}) \cdot \mathbf{v} \, dS &= 0
	\end{align}
	
	Due to the symmetry of strain rate  $\strainrate(\mathbf{u})$, we can write the inner product of strain rate tensor and test function $\mathbf{v}$ as:
	
	\[
	\strainrate(\mathbf{u}) : \nabla \mathbf{v} = \strainrate(\mathbf{u}) : \strainrate(\mathbf{v}).
	\]
	
	Thus, multiplying the incompressibility constraint by test function $q$ and integrating, the weak form problem is to find $(\mathbf{u}, p) \in \mathbf{V} \times Q$ such that
	\begin{align}
		\int_\Omega 2\mu \, \strainrate(\mathbf{u}) : \strainrate(\mathbf{v}) \, d\mathbf{x} - \int_\Omega p \, \nabla \cdot \mathbf{v} \, d\mathbf{x} - \int_\Omega \mathbf{f} \cdot \mathbf{v} \, d\mathbf{x} - \int_{\Gamma_N} (\boldsymbol{\sigma} \cdot \mathbf{n}) \cdot \mathbf{v} \, dS &= 0 \quad \forall \mathbf{v} \in \mathbf{V}_0 \label{eq:weak_momentum} \\
		\int_\Omega q \, \nabla \cdot \mathbf{u} \, d\mathbf{x} &= 0 \quad \forall q \in Q \label{eq:weak_incompressibility}
	\end{align}
	
	Jointly, Equations \eqref{eq:weak_momentum} and \eqref{eq:weak_incompressibility} are equivalent to those derived using a variational approach in \citet{boffiMixedFiniteElement2013}, with the exception of the boundary condition term in \eqref{eq:weak_momentum}. This term indicates that the total divergence of stress power over the domain must be equilibrated by the flux of stress power across the boundary, represented by the surface traction included as Neumann boundary conditions. The latter includes the possibility that $\boldsymbol{\sigma}  \mathbf{n} = \mathbf{0}$. 
	
	\section{Variational approach}
	
	Concepts from calculus of variation \cite{mitrCalculusVariationsGelfand2022, weinstockCalculusVariationsApplications1974} can be used to derive the weak form of the differential equations. The variational approach considers the domain in its entirety and uses the energy minimization to obtain a global description of the problem, which gives the generalized weak form of the Stokes equation. This is in contrast to the standard weak form of the Stokes equation, which is a local description of the flow at an infinitesimal point.  
	
	To obtain the weak form through the variational approach we apply the Gateaux derivative \cite{Gateaux1913}, i.e. the directional derivative, of the energy integral with respect to the state variables, e.g. velocity and pressure. Throughout this section, variations are taken with respect to admissible velocity and pressure fields belonging to the same function spaces as the primal solution and satisfying the prescribed boundary conditions. Variations are assumed to be sufficiently smooth to preserve the incompressibility constraint. We distinguish between variations of the fields within a fixed domain and variations induced by infinitesimal perturbations of the domain. In the following, only first-order contributions from the perturbation of the state variables are retained. Here we only consider the linear viscous problem for the sake of tractability and for comparison with previous approaches which also only considered linear viscosity \cite{kausStabilizationAlgorithmGeodynamic2010, roseStabilityAccuracyFree2017}.\\
	
	
	We start by defining the energy integral corresponding to incompressible  Stokes flow given by {Eq.(1.3.23)} in \cite{boffiMixedFiniteElement2013} as shown below:
	\begin{equation}
		I(\bm u, p) = \int_{\Omega(t)} F({\bf u}, \nabla {\bf u}, p) \,d\bold x = \int_{\Omega(t)}   \mu \, \strainrate(\mathbf{u}): \strainrate(\mathbf{u}) - \mathbf{f} \cdot \mathbf{u} - p \, (\nabla \cdot \mathbf{u})  \,d\bold x
		\label{minimization}
	\end{equation}
	
	here $\Omega(t)$ is the deformed domain as a function of time. Equation \eqref{minimization} represents the total internal energy balance of the system. Pressure $p$ also enters into the energy integral as a Lagrange multiplier to enforce the continuity equation.   The system is in equilibrium when the energy integral is minimized. Minimization is achieved by taking the derivative of the energy integral with respect to $\mathbf\alpha$, which controls the amount of perturbation away from the true solution as shown below:
	
	\begin{equation}
		\lim_{\alpha\to 0}\frac{I(\mathbf u + \alpha \mathbf ~\delta \mathbf u, p + \alpha~\delta p) - I(\mathbf u, p)}{\alpha} = 0
	\end{equation}
	
	We first present the solution for a fixed domain that does not include a free surface,
	and then the formulation for a deforming domain.

	\subsection{Variational approach for a fixed domain}
	
	In this section we consider the classical variational setting in which
	the \emph{domain is fixed} and only the state variables are perturbed.
	This contrasts with the next section, where both domain and field
	variations are present.
	Here, all variations are purely functional, and no geometric terms
	appear. We first define the general standard weak form of the Euler–Lagrange equation for
	functionals depending on first gradients, and then apply it to the Stokes equation.
	
	\medskip
	\noindent
	
	We consider a functional of the form
	\begin{equation}
		I(\mathbf{u})
		=
		\int_{\Omega}
		F\big(\mathbf{x}, \mathbf{u}(\mathbf{x}), \nabla \mathbf{u}(\mathbf{x})\big)\, d\Omega,
		\qquad
		\mathbf{u} : \Omega \to \mathbb{R}^d.
	\end{equation}
	The integrand $F$ is a function of state variable and its derivatives, namely the field
	$\mathbf{u}$, and its gradient.
	Since $\Omega$ is fixed, variations affect only $\mathbf{u}$ and
	$\nabla \mathbf{u}$ as it is usually defined in calculus of variations \cite{weinstockCalculusVariationsApplications1974, mitrCalculusVariationsGelfand2022}.
	
	\medskip
	\noindent
	
	Here we define the perturbation of a function $\mathbf{u}$ as:
	
	\begin{equation}
		\mathbf{u^\ast(\mathbf{x})} = \mathbf{u(\mathbf{x})} + \delta \mathbf{u(\mathbf{x})},
		\qquad
		\nabla \mathbf{u^\ast(\mathbf{x})}
		=
		\nabla \mathbf{u(\mathbf{x})} +  \delta \nabla( \mathbf{u(\mathbf{x})}),
		\label{perturbation_intro}
	\end{equation}
	
	where $\delta\mathbf{u}$ and $\delta \nabla {\mathbf{u}}$ are arbitrary admissible perturbations.
	
	\medskip
	\noindent
	
	The change in the functional induced by the perturbation is
	
	\begin{equation}
		\Delta I(\alpha)
		=
		\int_{\Omega}
		\Big[
		F\big(\mathbf{x}, \mathbf{u} + \alpha\,\delta \mathbf{u},
		\nabla \mathbf{u} + \alpha\,\nabla(\delta \mathbf{u})\big)
		-
		F(\mathbf{x},\mathbf{u},\nabla \mathbf{u})
		\Big]\,
		d\Omega.
		\label{deltaI}
	\end{equation}
	where $\alpha$ is a scalar. This expression measures the difference between the perturbed and
	unperturbed states.
	
	\medskip
	\noindent
	
	Expanding $F$ with respect to $\mathbf{u}$ and $\nabla\mathbf{u}$ to
	first order in $\alpha$ transforms \eqref{deltaI} into:
	\begin{equation}
		\Delta I(\alpha)
		=
		\int_{\Omega}
		\Big[
		\frac{\partial F}{\partial \mathbf{u}} \cdot (\alpha\,\delta \mathbf{u})
		+
		\frac{\partial F}{\partial (\nabla \mathbf{u})}
		:
		(\alpha\,\nabla(\delta \mathbf{u}))
		\Big]\,
		d\Omega.
	\end{equation}
	Higher-order terms in $\alpha$ are neglected, as we are interested in
	the linearized response of the functional.
	
	\medskip
	\noindent
	We minimize the action integral by using the directional derivative (Gâteaux derivative), which yields the first
	variation:
	\begin{equation}
		\delta I
		=
		\lim_{ \alpha \to 0} \frac{\Delta I}{\alpha}
		=
		\int_{\Omega}
		\left(
		\frac{\partial F}{\partial \mathbf{u}} \cdot \delta \mathbf{u}
		+
		\frac{\partial F}{\partial (\nabla \mathbf{u})}
		:
		\nabla(\delta \mathbf{u})
		\right)
		\, d\Omega.
		\label{general_functional}
	\end{equation}
	This is the standard weak form of the functionals depending on a function $\mathbf{u}$ and its gradients $\nabla \mathbf{u}$.
	
	\medskip
	\noindent
	\subsubsection{Application to Stokes equation.}
	We now apply the weak form of equation \eqref{general_functional}, to the energy integral of the incompressible Stokes flow, as defined by equation \eqref{minimization}, which yields:
	\medskip
	\noindent
	
	\begin{equation}
		\delta I
		=
		\int_{\Omega}
		\Big(
		2\mu\, \bm{\dot{\varepsilon}}(\mathbf{u})
		:
		\bm{\dot{\varepsilon}}(\delta \mathbf{u})
		- \mathbf{f} \cdot \delta \mathbf{u}
		- p\, \nabla \cdot (\delta \mathbf{u})
		- \delta p\, \nabla \cdot \mathbf{u}
		\Big)
		\, d\Omega.
	\end{equation}
	where the variations are taken with respect to the state variables $p$ and $\mathbf{u}$. The term involving $\delta p$ enforces the incompressibility constraint
	in a weak sense as an scalar field.
	
	\medskip
	\noindent
	
	We restrict variations to admissible test spaces, and renaming them accordingly as
	\begin{equation}
		\delta \mathbf{u} = \mathbf{v} ,
		\qquad
		\delta p = q
		\label{renaming}
	\end{equation}
	The homogeneous boundary condition on $\mathbf{v}$ corresponds to
	prescribed velocity on $\Gamma_D$.
	
	\medskip
	\noindent
	
	The stationary condition $\delta I=0$ for all admissible variations
	leads to
	\begin{equation}
		\int_{\Omega}
		\Big(
		2\mu\, \bm{\dot{\varepsilon}}(\mathbf{u}) : \bm{\dot{\varepsilon}}(\mathbf{v})
		- \mathbf{f} \cdot \mathbf{v}
		- p\, \nabla \cdot \mathbf{v}
		- q\, \nabla \cdot \mathbf{u}
		\Big)
		\, d\Omega
		= 0.
	\end{equation}
	This is the standard mixed variational formulation of the Stokes problem.
	
	\medskip
	\noindent
	
	Separating velocity and pressure test functions yields the momentum
	equation,
	\begin{equation}
		\int_{\Omega}
		\Big(
		2\mu\, \bm{\dot{\varepsilon}}(\mathbf{u}) : \bm{\dot{\varepsilon}}(\mathbf{v})
		- p\, \nabla \cdot \mathbf{v}
		- \mathbf{f} \cdot \mathbf{v}
		\Big)
		\, d\Omega
		= 0
		\qquad
		\forall \mathbf{v}.
	\end{equation}
	This represents the weak form of the balance of linear momentum.
	
	\medskip
	\noindent
	
	Finally, variation with respect to $p$ yields
	\begin{equation}
		\int_{\Omega}
		q\, \nabla \cdot \mathbf{u} \, d\Omega
		= 0
		\qquad
		\forall q,
		\label{p_variation}
	\end{equation}
	which enforces incompressibility in the weak sense.
	Together, the last two equations define the classical Stokes system on a
	fixed domain.

	\subsection{Variational approach for deforming domains}
\begin{figure}[h!]
		\centering
		\includegraphics[width=1\textwidth]{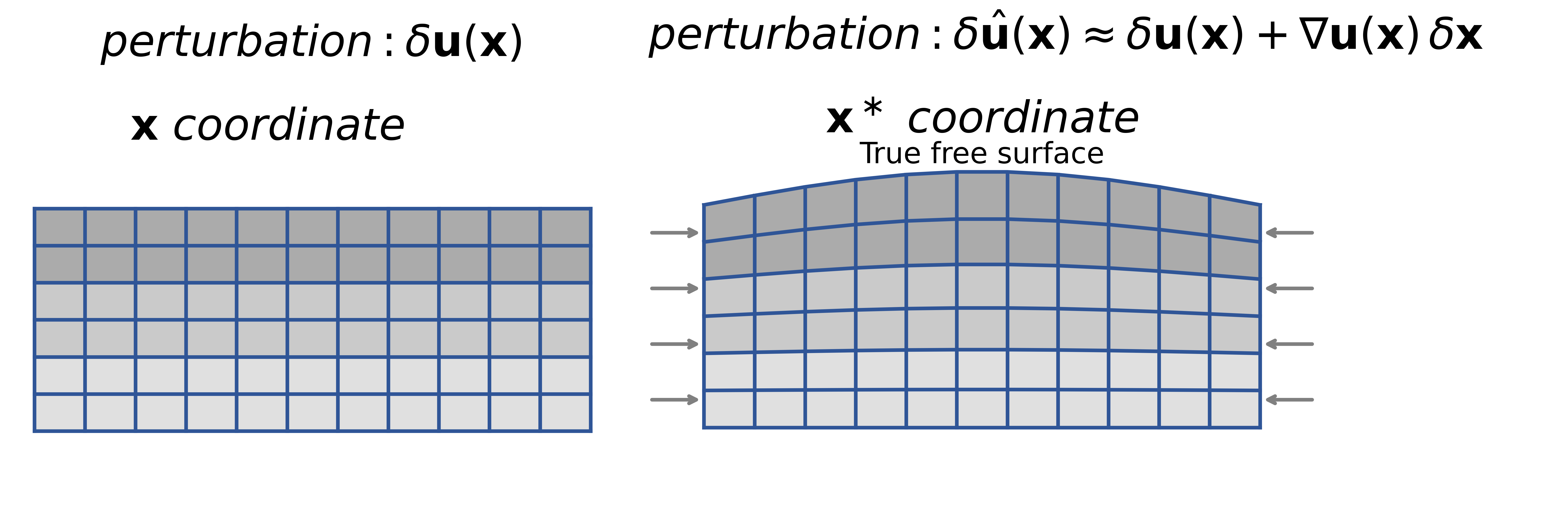}
		\caption{Free surface deformation modifies the domain geometry, requiring the state variables to be expressed in the deformed coordinate system. The problem is no longer steady state. }
		\label{fig:ale}
	\end{figure}

	We can perturb the domain and the state variables in order to produce the weak form for the free surface Stokes flow. In this setting, the future location of part of the domain associated with the free surface is also unknown. Although the Stokes flow does not have an explicit time dependence, free surface boundary allows for the deformation of domain, which in turn produces implicit time dependence.
	Therefore, we can use the domain location as a variable and apply a variation to it alongside the variation of the state variables. A schematic presentation of the free surface as well as the fixed and deformed coordinate systems are shown in Figure \ref{fig:ale}).
	
	We derive the first order variation of a energy functional when both the \emph{domain}
	and the \emph{state variables} are perturbed. The central idea is that a perturbation of the domain induces an additional variation to the state variables.
	Variations arising from perturbations of the domain are referred to as geometric variations. These are distinguished from variations caused by perturbations of the state variables that are independent of the domain perturbation; we refer to these as intrinsic variations.

	\medskip
	\noindent
	
	We start by defining an infinitesimal deformation of the spatial coordinates:
	\[
	\mathbf{x}^\ast = \mathbf{x} + \delta \mathbf{x}
	\label{deltax_star}
	\]
	Here $\mathbf{x}^\ast$ is the coordinate in the deformed domain and $\delta\mathbf{x}$ is a smooth field representing the domain displacement.
	This mapping sends the reference domain $\Omega$ to the deformed
	domain $\Omega_t$.
	
	\medskip
	\noindent
	
	The state variable is evaluated at different spatial points before and
	after the domain perturbation:
	
	The perturbed state variable in the deformed domain $\mathbf{u}^\ast(\mathbf{x}^\ast)$  is equal to the original value of the state variable, u(x), plus a variation that depends on both, geometric and intrinsic variations, therefore the perturbation in the deforming domain is distinguished by hat symbol on the state variable and $\delta \hat{\mathbf{u}}(\mathbf{x})$ is defined as:
	
	\begin{align}
		\mathbf{u}^\ast(\mathbf{x}^\ast)	
		&= \mathbf{u}(\mathbf{x}) +\delta \hat{\mathbf{u}}(\mathbf{x})
		\label{u_defrormed} \\
	\end{align}
	Therefore:
	\begin{align}
		\delta \hat{\mathbf{u}}(\mathbf{x})
		&= \mathbf{u}^\ast(\mathbf{x}^\ast) - \mathbf{u}(\mathbf{x})
		\label{u_star_perturbed} \\
		&= \mathbf{u}^\ast(\mathbf{x} + \delta \mathbf{x}) - \mathbf{u}(\mathbf{x})
		\label{u_star_deltax_perturbed}
	\end{align}
	the $\mathbf{u}^\ast(\mathbf{x}^\ast)$ is the perturbed state variable in the deformed domain. Equation \eqref{u_star_perturbed} defines the total variation of $\mathbf{u}$, while the second equation, 	\eqref{u_star_deltax_perturbed}, makes
	the change of coordinates explicit.
	
	To separate the intrinsic from the geometric variations, we expand
	$\mathbf{u}^\ast(\mathbf{x}+\delta\mathbf{x})$ in space.
	Using a first-order Taylor expansion we obtain:
	\begin{equation}
		\mathbf{u}^\ast(\mathbf{x} + \delta \mathbf{x})
		\approx \mathbf{u}^\ast(\mathbf{x})
		+ \nabla \mathbf{u}^\ast(\mathbf{x})\,\delta \mathbf{x}.
		\label{u_deform_var}
	\end{equation}
	We also define the intrinsic variation in the state variable $\mathbf{u}$ in the non-deforming frame as:  
	
	\begin{equation}
		\delta\mathbf{u}(\mathbf{x})=\mathbf{u}^\ast(\mathbf{x})-\mathbf{u}(\mathbf{x}) 
		\label{u_fixed_var}
	\end{equation}
	\\
	Using \eqref{u_fixed_var} and \eqref{u_deform_var} the total variation of $\hat{\mathbf{u}}(\mathbf{x})$ defined in \eqref{u_star_deltax_perturbed} becomes:
	
	\begin{equation}
		\delta \hat{\mathbf{u}}(\mathbf{x})
		\approx \delta \mathbf{u}(\mathbf{x})
		+ \nabla \mathbf{u}(\mathbf{x})\,\delta \mathbf{x},
		\label{eq:var_u}
	\end{equation}
	
	Equation~\eqref{eq:var_u} shows that the total variation consists of an
	intrinsic variation (first term), defined in the non-deforming frame, plus a geometric term (second term), induced by the domain deformation.

	\medskip
	\noindent
	
	The same reasoning applies to the gradient of the state variable.
	By definition,
	\begin{align}
		\nabla \delta \hat{\mathbf{u}}(\mathbf{x})
		&= \nabla \mathbf{u}^\ast(\mathbf{x}^\ast) - \nabla \mathbf{u}(\mathbf{x}) \\
		&= \nabla \mathbf{u}^\ast(\mathbf{x} + \delta \mathbf{x}) - \nabla \mathbf{u}(\mathbf{x}).
	\end{align}
	Expanding again to first order gives
	\begin{equation}
		\nabla \mathbf{u}^\ast(\mathbf{x} + \delta \mathbf{x})
		\approx \nabla \mathbf{u}^\ast(\mathbf{x})
		+ \nabla^2 \mathbf{u}^\ast(\mathbf{x})\,\delta \mathbf{x},
	\end{equation}
	and therefore
	\begin{equation}
		\nabla \delta\hat{\mathbf{u}}(\mathbf{x})
		\approx \delta \nabla \mathbf{u}(\mathbf{x})
		+ \nabla^2 \mathbf{u}(\mathbf{x})\,\delta \mathbf{x}.
		\label{eq:var_grad}
	\end{equation}
	Equation \eqref{eq:var_grad} shows that domain perturbations introduce
	second derivatives of the state variable, reflecting the geometric
	sensitivity of gradients. In the calculation of perturbation of the state variable and its gradient, we replaced the derivatives of $\mathbf{u}^\ast$ with derivatives of $\mathbf{u}$. This is acceptable as we will consider the limiting case when the directional derivative is used where the limiting case of $\alpha \to 0$ is considered. The reader may refer to \cite{mitrCalculusVariationsGelfand2022} for a more rigorous derivation of these formula.
	
	\medskip
	\noindent
	Now we consider the energy integral defined on the deformed domain,
	\begin{equation}
		I = \int_{\Omega_t} F(\mathbf{x}^\ast, \mathbf{u}^\ast, \nabla^\ast \mathbf{u}^\ast)\, d\Omega_t.
	\end{equation}
	To compute variations, it is convenient to pull this integral back to the
	fixed reference domain $\Omega$:
	\begin{equation}
		I = \int_{\Omega}
		F(\mathbf{x}^\ast, \mathbf{u}^\ast, \nabla \mathbf{u}^\ast)\, J\, d\Omega,
		\label{eq:I_integral}
	\end{equation}
	where the Jacobian of the transformation is an scalar value defined by \cite{malvernIntroductionMechanicsContinuous2007, jogContinuumMechanicsFoundations2015} as:
	\begin{equation}
		J = 1 + \alpha\,\nabla \cdot \delta \mathbf{x}.
	\end{equation}
	where $\alpha$ is the perturbation parameter  defined in \eqref{deltaI} and, $\nabla \cdot \delta \mathbf{x}$ is the divergence of $\delta \mathbf{x}$.
	This linearized Jacobian expresses the local change of volume induced by
	the domain deformation .
	
	\medskip
	\noindent
	
	Rewriting the energy integral equation \eqref{eq:I_integral} as a function of the perturbations in $\mathbf{x}$ and $\mathbf{u}$, we obtain:
	\begin{equation}
		I(\alpha)
		=
		\int_{\Omega}
		F\!\left(
		\mathbf{x} + \alpha\,\delta \mathbf{x},
		\mathbf{u}(\mathbf{x}) + \alpha \,\delta \hat{\mathbf{u}}(\mathbf{x}),
		\nabla \mathbf{u} + \alpha\,\delta \nabla \hat{\mathbf{u}}(\mathbf{x})
		\right)
		\bigl( 1 + \alpha\, \nabla \cdot \delta \mathbf{x} \bigr)
		\, d\Omega .
		\label{pullback}
	\end{equation}

	\medskip
	\noindent
	
	Expanding $F$ to first order using Taylor series expansion in $\alpha$ gives
	\begin{align}
		I(\alpha)
		&=
		\int_{\Omega}
		\Big[
		F(\mathbf{x},\mathbf{u},\nabla\mathbf{u})
		+ \alpha\, \frac{\partial F}{\partial \mathbf{x}} \cdot \delta \mathbf{x}
		+ \alpha\, \frac{\partial F}{\partial \mathbf{u}} \cdot \delta \hat{\mathbf{u}}(\mathbf{x})
		+ \alpha\, \frac{\partial F}{\partial (\nabla \mathbf{u})}
		:
		\nabla  \delta \hat{\mathbf{u}}(\mathbf{x})
		\Big]
		\bigl( 1 + \alpha\, \nabla \cdot \delta \mathbf{x} \bigr)
		\, d\Omega .
		\label{expand_f}
	\end{align}

	\medskip
	\noindent
	Using the relations derived in
	\eqref{eq:var_u}--\eqref{eq:var_grad}, we substitute
	\begin{align}
		\delta \hat{\mathbf{u}}(\mathbf{x})
		&= \delta \mathbf{u}(\mathbf{x}) + \nabla \mathbf{u} \, \delta \mathbf{x}
		\label{variations_total_u} \\[0.3em]
		\delta \nabla \hat{\mathbf{u}}(\mathbf{x})
		&=  \nabla \delta \mathbf{u}(\mathbf{x})
		+ \nabla^2 \mathbf{u}(\mathbf{x}) \, \delta \mathbf{x}
		\label{variations_total_grad_u}
	\end{align}
	\medskip
	\noindent
	
	After substitution of~\eqref{variations_total_u} and~\eqref{variations_total_grad_u} in~\eqref{expand_f} and collection of all first order ($\mathcal O (\alpha)$) contributions, we obtain
	\begin{align}
		I(\alpha)
		&=
		\int_{\Omega}
		\Big[
		F
		+ \alpha\, \frac{\partial F}{\partial \mathbf{x}} \cdot \delta \mathbf{x}
		+ \alpha\, \frac{\partial F}{\partial \mathbf{u}} \cdot \delta \mathbf{u}
		+ \alpha\, \frac{\partial F}{\partial \mathbf{u}} \cdot \nabla \mathbf{u}\,\delta \mathbf{x}
		\\
		&\qquad\quad
		+ \alpha\, \frac{\partial F}{\partial (\nabla \mathbf{u})}
		:
		(\nabla \delta \mathbf{u})
		+ \alpha\, \frac{\partial F}{\partial (\nabla \mathbf{u})}
		:
		\nabla^2 \mathbf{u}\,\delta \mathbf{x}
		+ \alpha\, F\, \nabla \cdot \delta \mathbf{x}
		\Big]
		\, d\Omega .
		\label{variation_allterms}
	\end{align}
	
	\medskip
	\noindent
	Because the divergence of $F\, \delta \mathbf{x}$ can be written as:
	
	\begin{equation}
		\nabla \cdot (F\,\delta \mathbf{x})
		=
		\frac{\partial F}{\partial \mathbf{x}} \cdot \delta \mathbf{x}
		+ \frac{\partial F}{\partial \mathbf{u}} \cdot (\nabla \mathbf{u}\,\delta \mathbf{x})
		+ \frac{\partial F}{\partial \nabla \mathbf{u}} : (\nabla^2 \mathbf{u}\,\delta \mathbf{x})
		+ F\,\nabla \cdot \delta \mathbf{x}
	\end{equation}
	all terms proportional to $\delta\mathbf{x}$ combine into a divergence in equation \eqref{variation_allterms} as:
	\begin{align}
		I(\alpha) 
		&=
		\int_{\Omega}
		\Big[
		F
		+ \alpha\, \nabla \cdot ( F\, \delta \mathbf{x} )
		+ \alpha\, \frac{\partial F}{\partial \mathbf{u}} \cdot \delta \mathbf{u}
		+ \alpha\, \frac{\partial F}{\partial (\nabla \mathbf{u})}
		:
		\delta (\nabla \mathbf{u})
		\Big]
		\, d\Omega .
	\end{align}

	which shows that the \textbf{geometric variation reduces to a divergence}.
	\medskip
	\noindent
	
	Taking the derivative with respect to $\alpha$ at $\alpha=0$ yields
	\begin{equation}
		\delta I
		=
		\int_{\Omega}
		\left(
		\frac{\partial F}{\partial \mathbf{u}} \cdot \delta \mathbf{u}
		+
		\frac{\partial F}{\partial (\nabla \mathbf{u})}
		:
		\delta (\nabla \mathbf{u})
		\right)
		\, d\Omega
		+
		\int_{\Omega}
		\nabla \cdot ( F\, \delta \mathbf{x} ) \, d\Omega
		\label{canonical}
	\end{equation}
	The first integral corresponds to the classical variational derivative
	with respect to the state variable, while the second integral captures the geometric
	perturbation.
	
	\medskip
	\noindent
	
	Applying the divergence theorem once more yields
	\begin{equation}
		\int_\Omega \nabla \cdot (F\,\delta \mathbf{x})\, d\Omega
		= \int_{\Gamma_N} F\,\delta \mathbf{x}\cdot \mathbf{n}\, dS
		\label{divergence}
	\end{equation}
	where $\mathbf{n}$ the normal vector at the boundary of domain.
	Application of identity \eqref{divergence} onto equation \eqref{canonical} produces the final result as:
	\begin{equation}
		\delta I
		=
		\int_{\Omega}
		\left(
		\frac{\partial F}{\partial \mathbf{u}} \cdot \delta \mathbf{u}
		+
		\frac{\partial F}{\partial (\nabla \mathbf{u})}
		:
		\delta (\nabla \mathbf{u})
		\right)
		\, d\Omega
		+
		\int_{\Gamma_N}
		F\, \delta \mathbf{x} \cdot \mathbf{n}
		\, dS
		\label{eq:result_canonical}
	\end{equation}
	
	Therefore, domain variations contribute only through boundary integral. The result explicitly displays both intrinsic  and geometric contributions. It is worth noting that the  equation \eqref{eq:result_canonical} reduces to the fixed domain variation when the domain perturbation vanishes.
	
	Equation \eqref{eq:result_canonical} is the standard result from shape sensitivity analysis \cite{sokolowskiIntroductionShapeOptimization2012}, showing that domain variations contribute only through a boundary integral. For a comprehensive treatment of the geometric aspects of shape variation, we refer to Henrot and Pierre \cite{henrotShapeVariationOptimization2018}.
	We will apply this equation in a novel way to the minimization of energy integral of the Stokes flow with free surface in the next section.
	
	\subsubsection{Application to Stokes equation with free surface.}
	
	In a similar fashion to what we did to the fixed domain, we apply the state variable and spatial perturbation to the energy integral equation \eqref{minimization}.

	
	\begin{align}
		& \int_{\Omega}
		\left(
		2 \mu \, \strainrate_{\mathbf{u}} : \strainrate_{\mathbf{\delta{u}}}
		- \mathbf{f} \cdot \mathbf{\delta{u}}
		- p \, \nabla \cdot \mathbf{\delta{u}}
		- \delta p \, \nabla \cdot \mathbf{u}
		\right)
		\, d\Omega
		\notag \\
		&\quad
		+
		\int_{\Gamma_N}
		\left(
		\mu \, \strainrate_{\mathbf{u}} : \strainrate_{\mathbf{u}}
		- \mathbf{f} \cdot \mathbf{u}
		- p \nabla \cdot \mathbf{u}
		\right)
		(\mathbf{\delta{\mathbf{x}}} \cdot \mathbf{n})
		\, dS
		= 0 .
		\label{eq:final}
	\end{align}

	The weak form in equation \eqref{eq:final} contains three types of perturbations, namely the intrinsic perturbation of the velocity and pressure field, and the geometric perturbation of the spatial domain.
	The perturbation of the state variables and the geometric perturbation are arbitrary and independent in general. Variation with respect to the pressure yields
	\[
	\int_{\Omega} \delta p \, \nabla \cdot \mathbf{u} \, d\Omega = 0
	\quad \forall \, \delta p,
	\]
	which implies the conservation of mass constraint:
	\[
	\nabla \cdot \mathbf{u} = 0 \quad \text{in } \Omega.
	\]
	Consequently, the term $p \nabla \cdot \mathbf{u}$ is identically zero, and in particular does not contribute to the boundary integral.
	
	Now, we can choose to relate the geometric perturbation of the domain to a perturbation of the velocity  by choosing the geometric perturbation to be:
	
	\begin{align}
		\delta{\mathbf{x}} = \Delta t \, \delta{\mathbf{u}}
		\label{eq:spatial_perturbation_reduction}
	\end{align}
	
	Using equation \eqref{eq:spatial_perturbation_reduction}, we reduce the perturbation space of the deforming domain to that of the fixed domain. Using equation \eqref{eq:spatial_perturbation_reduction} we can simplify the boundary integral as the term $p \nabla \cdot \mathbf{u}$ now vanishes. The physical interpretation is that the domain can undergo deformation in such a way that velocity remains divergence free.
	\begin{subequations}
		\begin{align}
			\int_{\Omega}
			\left(
			2 \mu \, \strainrate_{\mathbf{u}} : \strainrate_{\delta \mathbf{u}}
			- \mathbf{f} \cdot \delta \mathbf{u}
			- \delta p \, \nabla \cdot  \mathbf{u}
			\right)
			\, \mathrm{d}\Omega
			\notag \\
			\quad
			+ \Delta t
			\int_{\Gamma_N}
			\left(
			\mu \, \strainrate_{\mathbf{u}} : \strainrate_{\mathbf{u}}
			- \mathbf{f} \cdot \mathbf{u}
			\right)
			(\delta \mathbf{u} \cdot \mathbf{n})
			\, \mathrm{d}S
			= 0
			\label{canonical2a} \\
			\int_{\Omega}
			\left(
			- \delta p \, \nabla \cdot \mathbf{u}
			\right)
			\, \mathrm{d}\Omega
			= 0 .
			\label{canonical2b}
		\end{align}
	\end{subequations}

	
	For nearly horizontal free surfaces, the term $(\mathbf{f} \cdot \mathbf{u})$ may be replaced by $(\rho g \, \mathbf{n} \cdot \mathbf{u})$, yielding a symmetric boundary contribution in the weak form.
	
	The boundary integral in equation \eqref{canonical2b} contains two terms arising from the perturbation of the domain. The gravity term $\mathbf{f \cdot u}$ is commonly used and interpreted as a stabilization of free surface \cite{kausStabilizationAlgorithmGeodynamic2010, roseStabilityAccuracyFree2017}, while the \textbf{viscous term was never considered}. If viscous effects at the boundary are negligible, the latter term may be omitted, recovering the stabilization commonly employed in free-surface Stokes flow. In section 5, we explore where viscous effects at the boundary are negligible or not.
	
	In the limit when $\Delta t \to 0$, the boundary term vanishes and the formulation reduces to the standard steady state weak form of the Stokes equations. For larger time steps, the boundary contribution becomes non negligible and accounts for finite free surface deformation. Hence, the formulation is consistent in both the small and large time steps $\Delta t$.
	Using the same change of test function names as for the fixed domain according to equation \eqref{renaming}, we can rewrite the weak form as:
	\begin{subequations}
		\begin{align}
			\int_{\Omega}
			\left(
			2 \mu \, \strainrate_{\mathbf{u}} : \strainrate_{\mathbf{v}}
			- \mathbf{f} \cdot \mathbf{v}
			- p \, \nabla \cdot \mathbf{v}
			\right)
			\, \mathrm{d}\Omega
			\notag \\
			\quad
			+ \Delta t
			\int_{\Gamma_N}
			\left(
			\mu \, \strainrate_{\mathbf{u}} : \strainrate_{\mathbf{u}}
			- \mathbf{f} \cdot \mathbf{u}
			\right)
			(\mathbf{v} \cdot \mathbf{n})
			\, \mathrm{d}S
			= 0 ,
			\label{canonical2avq} \\
			\int_{\Omega}
			\left(
			- q \, \nabla \cdot \mathbf{u}
			\right)
			\, \mathrm{d}\Omega
			= 0 .
			\label{canonical2bvq}
		\end{align}
	\end{subequations}
	
	
	The boundary integral is approximated using the trapezoidal rule, corresponding to an average over the initial and final states of the deformation interval $\Delta t$. We note that While \citet{kausStabilizationAlgorithmGeodynamic2010} derived the gravity term through heuristic reasoning, our derivation shows it emerges naturally from the domain variation. Furthermore, we identify and quantify the previously overlooked viscous boundary term. The Firedrake framework \cite{rathgeberFiredrakeAutomatingFinite2017} is used to implement the solver for the above weak from. The details of the implementation are shown in the supplementary material. 
	
	
	
	\section{Dimensional analysis of the significance of terms}
	
	In order to assess the relative importance of the terms in the boundary integral for geodynamics problems, we perform a dimensional analysis. Dimensionless variables and parameters are denoted with $^*$, and the notation is restricted to this section.
	
	We start with the original functional:
	\begin{equation}
		I = \int_{\Omega} \left( \mu (\nabla \mathbf{u} + \nabla \mathbf{u}^T) : (\nabla \mathbf{u} + \nabla \mathbf{u}^T) - \rho \mathbf{g} \cdot \mathbf{u} \right) d\mathbf{x}.
	\end{equation}
	
	Introducing the following dimensionless variables:
	\begin{equation}
		\mathbf{x}^* = \frac{\mathbf{x}}{L}, \quad
		\mathbf{u}^* = \frac{\mathbf{u}}{U}, \quad
		\nabla^* = L \nabla, \quad
		\mathbf{g}^* = \frac{\mathbf{g}}{g},
	\end{equation}
	so that
	\begin{equation}
		\mathbf{u} = U \mathbf{u}^*, \quad
		\nabla = \frac{1}{L} \nabla^*, \quad
		d\mathbf{x} = L^d d\mathbf{x}^*.
	\end{equation}
	
	Using these relations, the viscous and gravity terms transform as
	\begin{align}
		(\nabla \mathbf{u} + \nabla \mathbf{u}^T) : (\nabla \mathbf{u} + \nabla \mathbf{u}^T)
		&= \frac{U^2}{L^2} (\nabla^* \mathbf{u}^* + \nabla^{*} \mathbf{u}^{*T}) : (\nabla^* \mathbf{u}^* + \nabla^{*} \mathbf{u}^{*T}), \\
		\rho \mathbf{g} \cdot \mathbf{u} &= \rho g U \mathbf{g}^* \cdot \mathbf{u}^*.
	\end{align}
	
	Substituting into the integral and factoring out $\rho U^3 L^d$ and $ \frac{1}{L}$ gives:
	\begin{equation}
		I = \rho U^3 L^{d-1} \int_{\Omega^*} \left( \frac{\mu}{\rho U L} (\nabla^* \mathbf{u}^* + \nabla^{*} \mathbf{u}^{*T}) : (\nabla^* \mathbf{u}^* + \nabla^{*} \mathbf{u}^{*T}) - \frac{g L}{U^2} \mathbf{g}^* \cdot \mathbf{u}^* \right) d\mathbf{x}^*.
	\end{equation}
	
	Introduce the Reynolds and Froude numbers:
	\begin{equation}
		Re = \frac{\rho U L}{\mu}, \quad Fr = \frac{U}{\sqrt{g L}},
	\end{equation}
	which allows rewriting the prefactor as
	\begin{equation}
		\frac{\mu}{\rho U L} = \frac{1}{Re}, \quad
		\frac{g L}{U^2} = \frac{1}{Fr^2}.
	\end{equation}
	
	Finally, we define the purely dimensionless functional $I^*$:
	\begin{equation}
		I^* = \int_{\Omega^*} \left( \frac{1}{Re} (\nabla^* \mathbf{u}^* + \nabla^{*} \mathbf{u}^{*T}) : (\nabla^* \mathbf{u}^* + \nabla^{*} \mathbf{u}^{*T}) - \frac{1}{Fr^2} \mathbf{g}^* \cdot \mathbf{u}^* \right) d\mathbf{x}^*,
	\end{equation}
	so that the original functional can be expressed as
	\begin{equation}
		I = \rho U^3 L^{d-1} \, I^*.
	\end{equation}
	
	This formulation makes explicit the scaling of the functional with $\rho$, $U$, and $L$, while $I^*$ is dimensionless. Since the same terms appear in the boundary integral of the weak form, this allows a quantitative estimate of the relative contributions of viscous and gravitational forces in various geodynamic problems.

	\begin{table}[htbp]
		\centering
		\caption{Characteristic scales for representative geodynamic systems. The ratio $Re/Fr^2 = \rho g L^2 / (\mu U)$ quantifies the relative importance of gravitational forcing to viscous resistance in the Stokes regime. Values in the final column are rounded to the nearest order of magnitude.}
		\label{tab:reduced_scaling}
		\begin{tabular}{lccccc}
			\hline
			\makecell{\textbf{Geodynamic} \\  \textbf{Problem}} &
			\makecell{\textbf{Length}\\(m)} &
			\makecell{\textbf{Velocity}\\(m/s)} &
			\makecell{\textbf{Density}\\(kg/m$^3$)} &
			\makecell{\textbf{Viscosity}\\(Pa$\cdot$s)} &
			\makecell{$\mathbf{Re/Fr^2}$} \\
		\hline
		Mantle Convection    & \(1.0\times10^{6}\)   & \(1\times10^{-9}\)  & \(3.30\times10^{3}\)  & \(1.0\times10^{21}\)  & \(3\times10^{4}\)  \\
		Rifted Margins       & \(1.0\times10^{5}\)   & \(1\times10^{-10}\) & \(3.30\times10^{3}\)  & \(1.0\times10^{22}\)  & \(3\times10^{2}\)  \\
		Ice Sheet Flow       & \(1.0\times10^{5}\)   & \(3\times10^{-8}\)  & \(9.00\times10^{2}\)  & \(1.0\times10^{14}\)  & \(3\times10^{7}\)  \\
		Plate Subduction     & \(5.0\times10^{5}\)   & \(3\times10^{-9}\)  & \(3.30\times10^{3}\)  & \(1.0\times10^{22}\)  & \(3\times10^{2}\)  \\
		\hline
		\end{tabular}
	\end{table}

	Table \ref{tab:reduced_scaling} shows some of geodynamics flows and the ranges of the expected ranges of Reynolds and Froude numbers determine the importance of the viscosity term and the gravity term in the boundary integral. It can be seen that the viscous term is important for Rifting and Plate subduction, and also should not be neglected in mantle convection and ice sheet and glacial systems. For systems with lower viscosity, such as lava flow and magma ascent, the viscous term can be ignored. For systems with smaller viscosity such as lava flow and magma ascent, can be ignored. \\
	It is worth noting that dimensional analysis provides a useful first-order estimate of the relative importance of the terms in the governing equation. Consequently, dimensional analysis should be regarded as a preliminary assessment rather than a definitive criterion. In particular, it can be used to identify and exclude limiting cases in which viscous effects are clearly negligible, but it may not reliably capture their influence on the overall model evolution.

	\section{Numerical Test Cases}
	\label{sec:test-cases}
	
	Two test cases are considered to assess the accuracy, stability, and
	applicability of the proposed free-surface stabilization formulations.
	The tests are arranged in increasing order of physical and numerical
	complexity. The first case is the Rayleigh--Taylor benchmark  and the second test is the Rising Plume test.
	
	\subsection{Rayleigh--Taylor Instability}
	\label{subsec:rayleigh-taylor-test}
	
\begin{figure}[h!]
		\centering
		\includegraphics[width=0.5\textwidth]{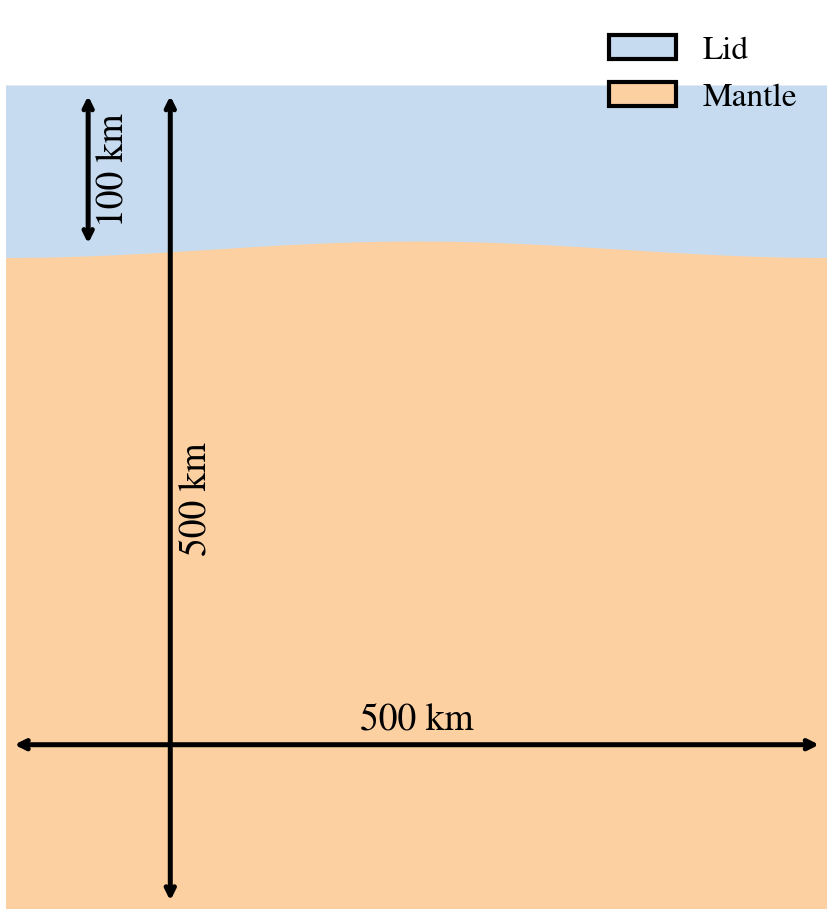}
		\caption{
			Rayleigh--Taylor instability benchmark. A denser and more viscous
			layer overlies a less dense and less viscous layer, producing an
			gravitationally unstable configuration. The test is used to assess
			the stability and accuracy of the free-surface formulation under
			increasing time-step sizes.
		}
		\label{fig:benchmark2}
	\end{figure}

	The first test case is the Rayleigh--Taylor instability, in which a
	denser, more viscous layer overlies a less dense and less viscous layer.
	The configuration represents an unstable gravitational arrangement and
	provides a standard test of the stability and accuracy of free-surface
	methods. The initial configuration is shown in
	Figure~\ref{fig:benchmark2}. The problem has been extensively used for
	the assessment of numerical methods for geodynamic free-surface
	calculations \cite{kausStabilizationAlgorithmGeodynamic2010,
		roseStabilityAccuracyFree2017}.
	
	The primary objective of this test is to examine whether the proposed
	stabilization terms suppress the numerical oscillations associated with
	the free surface while preserving the physical growth of the instability.
	The test is particularly sensitive to the time-step size: sufficiently
	large time steps may result in oscillatory free-surface motion and
	eventual divergence from the physically evolving solution. A series of
	simulations with different viscosity contrasts and time-step sizes is
	therefore used to evaluate the stability range of the competing
	formulations.
	
	\subsection{Rising Plume}
	\label{subsec:rising-plume-test}
	
	The second test case considers a buoyant plume rising through a more
	viscous surrounding medium and interacting with a deformable free
	surface. The configuration follows the idealized plume benchmark used
	by \citet{crameriComparisonNumericalSurface2012} to investigate errors
	in predicted surface topography associated with the sticky-air
	approximation.
	
	In contrast to that approach, the present simulations employ a true free
	surface and track its motion in a Lagrangian framework. The benchmark
	therefore provides a direct assessment of the ability of the proposed
	stabilization schemes to reproduce the surface deformation generated by
	an evolving buoyant structure without introducing an artificial
	low-viscosity layer.
	
	\begin{figure}[h!]
		\centering
		\includegraphics[width=\textwidth]{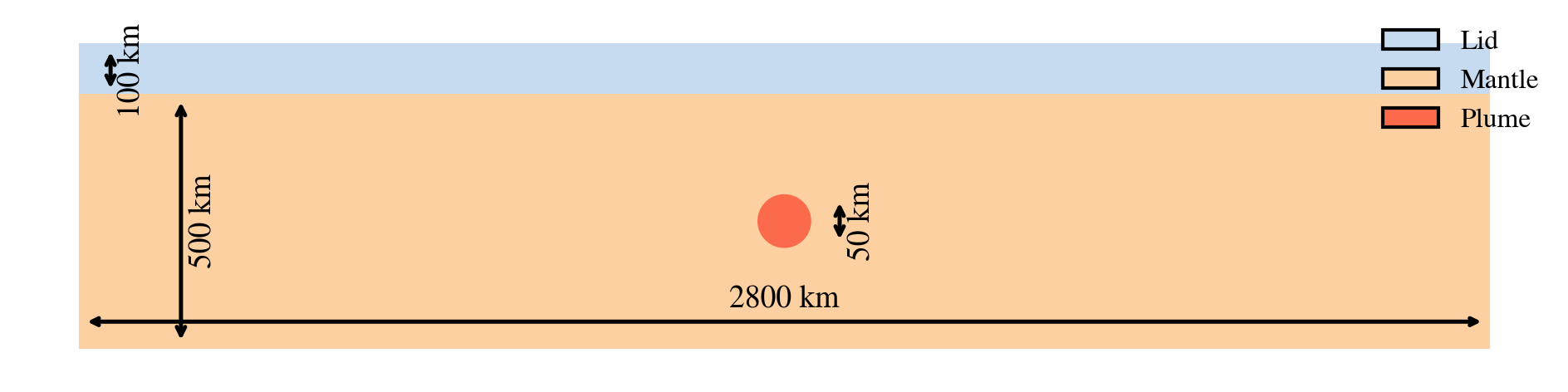}
		\caption{
			Rising-plume benchmark. A buoyant viscous inclusion rises through
			the surrounding material and deforms the free surface. The
			resulting surface topography depends on the relative viscosity of
			the plume and surrounding layers and provides a sensitive measure
			of the performance of the free-surface formulation.
		}
		\label{fig:rp-setup}
	\end{figure}
	
	The benchmark is used to compare the predicted surface evolution and
	topographic response of the different stabilization strategies. Because
	the viscosity is constant within each material in this test, the role of
	the stabilization terms can be assessed without the additional
	complications associated with a strongly evolving nonlinear rheology.
	The lid viscosity is varied systematically to investigate the robustness
	of the formulations over a range of viscosity contrasts.
	\pagebreak

	\section{Results}
	
	Stabilization of the free surface using the gravity term has been derived and studied heuristically in geodynamics literature \cite{crameriComparisonNumericalSurface2012, kausStabilizationAlgorithmGeodynamic2010, andres-martinezNewFreesurfaceStabilization2015}. Here we show that the Full stabilizationilization due to the variational approach to the Stokes problem is consistent with the already existing stabilization method for the problems studied as benchmarks in the free‑surface literature. In the test study cases, we show that the solutions are close to each other and the fully stabilized formulation is consistent with the previously used methods, validating the approach.
	
\begin{figure}[h!]
		\centering
		\includegraphics[width=\textwidth]{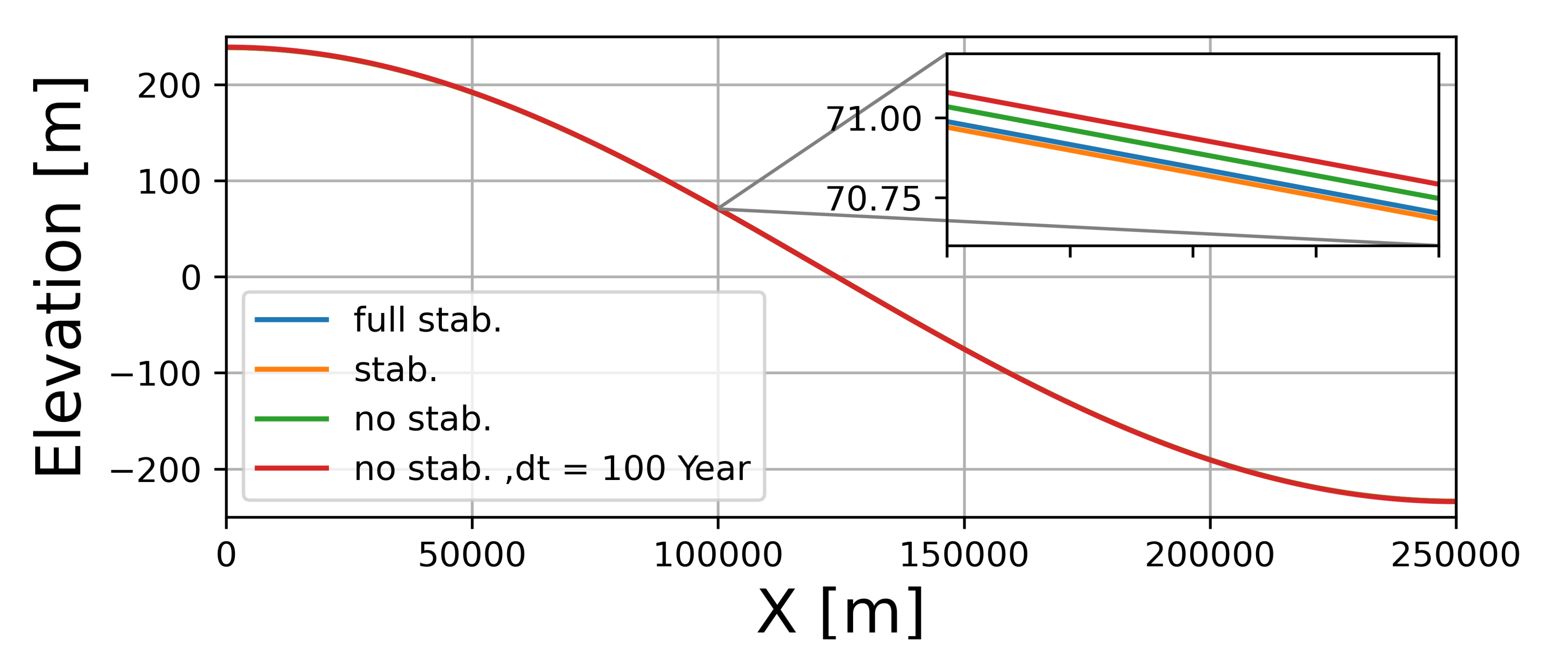}
		\caption{Evolution of the top interface in the Rayleigh–Taylor instability test case. The interface position is shown after 2~million years of simulation time. Comparison is made between the proposed Full stabilizationilization (incorporating viscous and gravity terms), gravity‑only stabilization, and a reference simulation, all using a time step of $dt=5000$~years. The reference solution with $dt=100$~years is provided for accuracy assessment. Results demonstrate that both stabilization techniques maintain high fidelity to the reference model within standard error ranges.}
		\label{fig:rt_topo}
	\end{figure}

\begin{figure}[H]
		\centering
		\includegraphics[width=\textwidth]{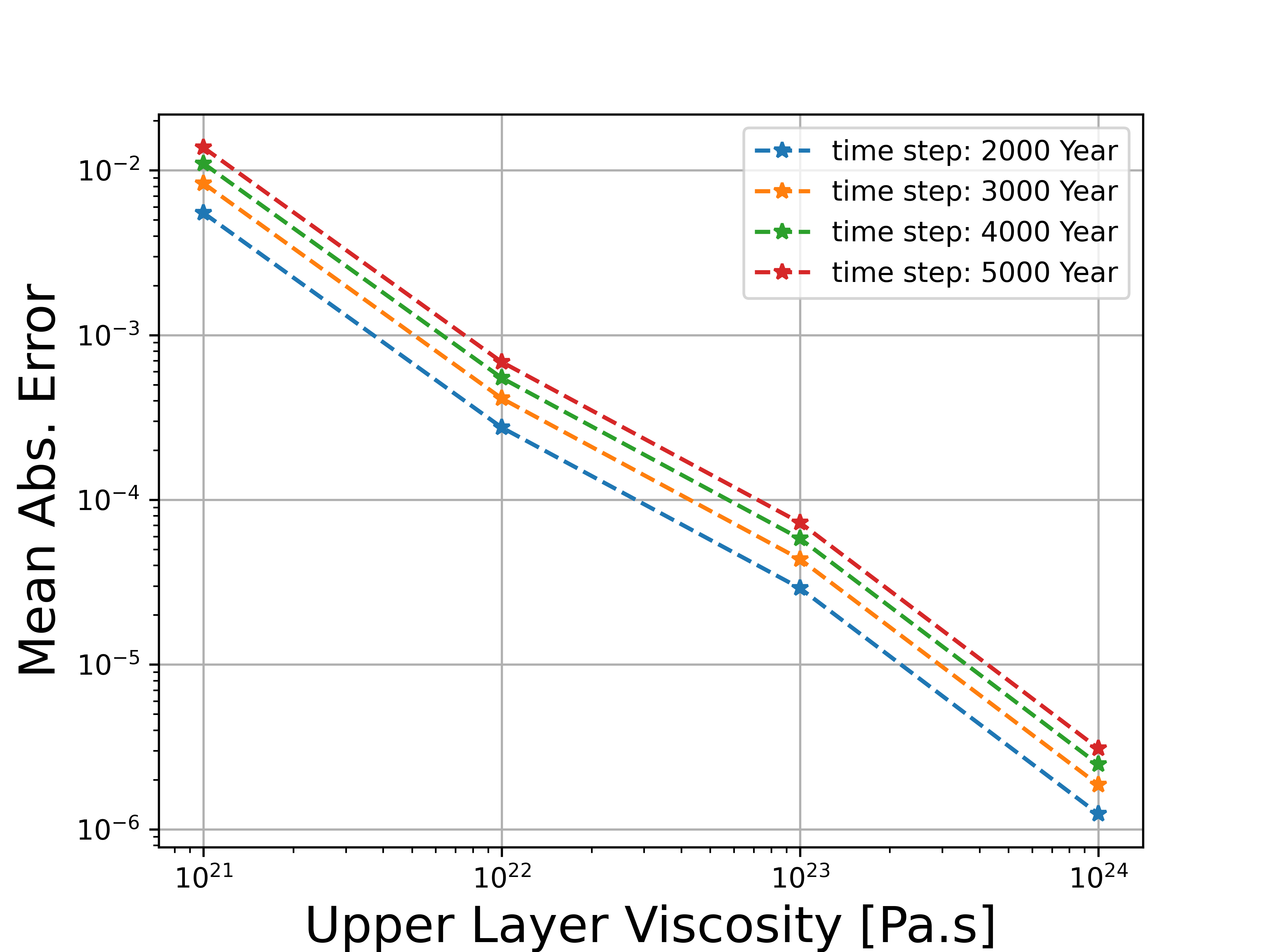}
		\caption{Mean Absolute Error (MAE), defined by Equation~\eqref{eq:mae}, for the upper surface of the Rayleigh–Taylor instability. The MAE is compared for both full and gravity‑only stabilization methods across various time steps and lid viscosities. Both methods demonstrate consistent performance. As expected, the MAE is inversely proportional to lid rigidity, peaking at lower viscosity values where surface deformation is most significant.}
		\label{fig:rt_meanerror}
	\end{figure}

	For the first benchmark we look at different time steps $[2000, 3000, 4000 ,\, 5000]$~years. The resulting surface topography is shown in Figure~\ref{fig:rt_topo}. The mean absolute error, defined by equation~\eqref{eq:mae} and graphically shown in Figure~\ref{fig:rt_meanerror}, is presented for various Reynolds numbers (different lid viscosity values).
	
	\begin{equation}
		\text{Mean Abs. Error} = \frac{1}{n} \sum_{i=1}^{n} |y_i - y_{oi}|
		\label{eq:mae}
	\end{equation}
	
	where $y_i - y_{oi}$ is the difference between large $\Delta t$ models and the reference model with $dt = 100$~years, and $n$ is the number of mesh nodes located at the surface boundary of the reference model.
	
	The deviation of surface topography from a reference model (calculated with a smaller time step of $dt = 100$~years) is determined using Equation~\eqref{eq:mae}. Figure~\ref{fig:rt_topo} illustrates the resulting relief for a simulation with a lid viscosity of $\mu = 10^{21}$~Pa$\cdot$s and a time step of $dt = 5000$~years. The results indicate that the Full stabilizationilization method and the gravity‑only stabilization method yield nearly identical results.
	
	We also examine the impact of varying lid viscosities on the solution. Figure~\ref{fig:rt_meanerror} displays the mean error between the fully stabilized and gravity‑stabilized variants, alongside a reference model with $dt = 100$~years. For the time steps evaluated, both stabilization methods perform consistently across all lid viscosities. 
	
	As expected, increasing the lid viscosity reduces overall deformation, causing both methods to converge towards a similar free surface. While larger discrepancies occur at lower viscosities (particularly $\mu = 10^{21}$~Pa$\cdot$s), these differences remain significantly smaller than the typical error margins for these model types. Furthermore, the observed variance is substantially lower than that reported for sticky‑air methods \cite{crameriComparisonNumericalSurface2012}.
	
\begin{figure}[h!]
		\centering
		\includegraphics[width=\textwidth]{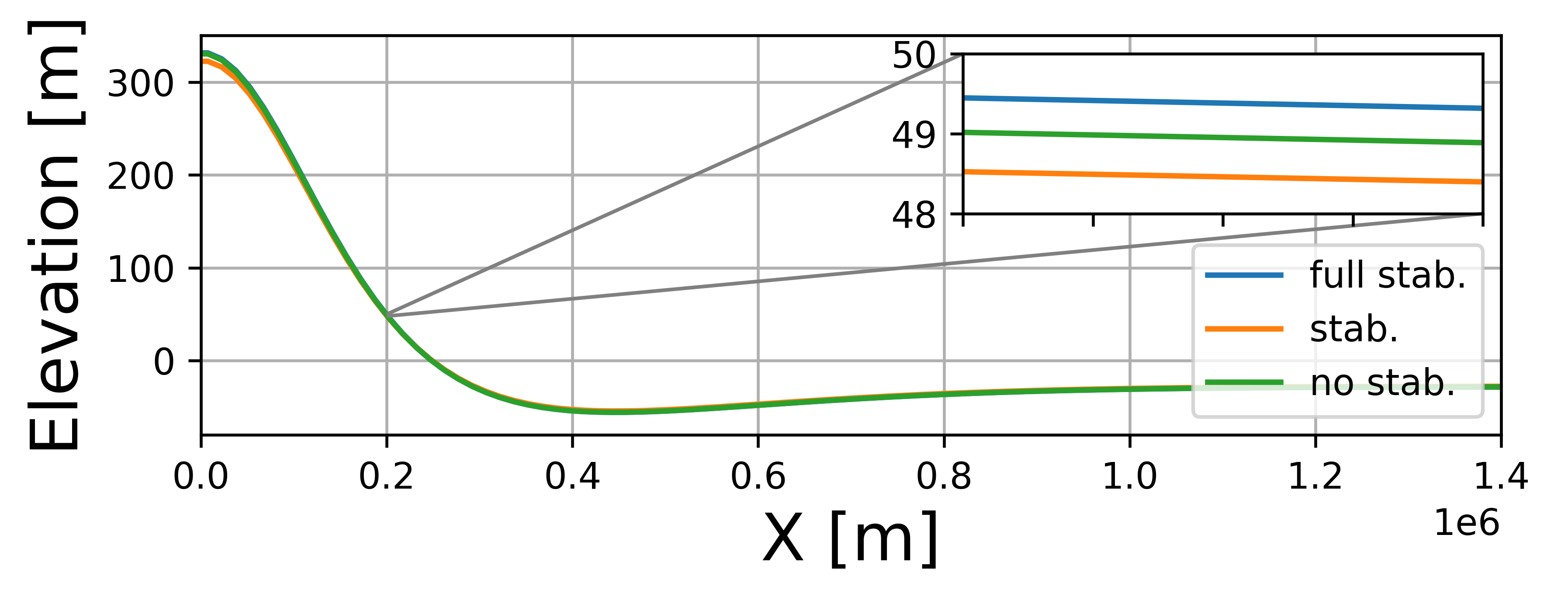}
		\caption{The upper surface is shown for half of the domain for the rising plume test case. The rising plume demonstrates a simple case of upward movement of lighter materials in the mantle. The free surface deforms accordingly. Here a fully stabilized, stabilized, and non‑stabilized model are shown. The inset shows a smaller part of the graph. Lid viscosity is $\mu = 1e21$ and $\Delta t = 5000$~years.}
		\label{fig:rp_topo}
	\end{figure}

	In the second case study, we examine the resulting topography for a lid viscosity of $\mu = 10^{21}$~Pa$\cdot$s and a time step of $dt = 5000$~years. Figure~\ref{fig:rp_topo} compares the fully stabilised, stabilised, and unstabilised simulations of the rising plume benchmark. Both stabilization methods yield results closely aligned with each other and accurately track the reference simulation ($dt = 100$~years).
	
	We also analyse the effect of varying the lid viscosity. For each viscosity value, the three simulation types remain relatively consistent. Table~\ref{tab:risingplume} summaries the maximum elevation, minimum elevation, and total relief after three million years. As lid viscosity increases (and the Reynolds number decreases), the three methods demonstrate improved agreement regarding maximum elevation. Conversely, the minimum elevation exhibits less variability across all simulations. The total relief peaks at a viscosity of $\mu = 10^{23}$~Pa$\cdot$s, while remaining comparable for $\mu = 10^{22}$ and $10^{24}$~Pa$\cdot$s. These results confirm that neither stabilization method introduces significant errors into the predicted topography.
	
	For both test cases, the label \q{Full stabilization} denotes the new derivation incorporating both viscous and gravity effects, while \q{Gravity-only stabilization} refers to the application of the gravity term alone within the boundary integral.
	
	\begin{table}[ht!]
		\centering
		\begin{tabular}{|c|c|c|c|c|}
			\hline
			\textbf{stabilization} 
			& \makecell{\textbf{Viscosity} \\ (Pa $\cdot$ s)} 
			& \makecell{\textbf{Max Elevation} \\ (km)} 
			& \makecell{\textbf{Min Elevation} \\ (km)} 
			& \makecell{\textbf{Relief} \\ (km)} \\
			\hline
			Full stabilization & $10^{21}$ & 224.31 & -38.50 & 262.81 \\
			stabilization      & $10^{21}$ & 226.62 & -38.81 & 265.43 \\
			No stabilization   & $10^{21}$ & 208.16 & -35.79 & 243.95 \\
			\hline\hline
			Full stabilization & $10^{22}$ & 286.16 & -43.95 & 330.11 \\
			stabilization      & $10^{22}$ & 309.35 & -47.72 & 357.07 \\
			No stabilization   & $10^{22}$ & 302.30 & -46.79 & 349.09 \\
			\hline\hline
			Full stabilization & $10^{23}$ & 314.10 & -43.47 & 357.57 \\
			stabilization      & $10^{23}$ & 327.86 & -45.23 & 373.09 \\
			No stabilization   & $10^{23}$ & 300.75 & -40.35 & 341.10 \\
			\hline\hline
			Full stabilization & $10^{24}$ & 270.11 & -60.47 & 330.58 \\
			stabilization      & $10^{24}$ & 267.91 & -60.61 & 328.52 \\
			No stabilization   & $10^{24}$ & 269.36 & -60.81 & 330.18 \\
			\hline
		\end{tabular}
		\caption{Relief statistics for full, gravity‑only, and unestablished formulations for the rising plume case after $3$~Myr, using a time step of $dt = 5000$~years. Both full and gravity‑only stabilization yield consistent results for minimum elevation, maximum elevation, and total relief across the range of viscosities considered.}
		\label{tab:risingplume}
	\end{table}
	

	\section{Discussion and Conclusion}
	
	We showed that the boundary integral on the free surface commonly used as the stabilization term for the Stokes equations is a natural product of the perturbation of the domain. The perturbation of the domain reveals that there is another stabilization term due to viscosity. Dimensional analysis shows that the importance of the gravity and viscosity terms are proportional to the Froude and Reynolds numbers, respectively. The analysis also shows that the viscosity term can be significant in many problems. For the simple test cases studied in the literature, which are also used in this study, the viscosity term has a smaller significance.
	
	It should be emphasized that the present formulation was derived under the assumption of a linear‑viscosity problem. 
	
	A comprehensive empirical assessment of the viscous stabilization term remains beyond the current study's scope. The dynamics of a free surface may show high sensitivity to the chosen time‑integration scheme (e.g., Runge–Kutta or Adams–Bashforth) for the boundary integral. Consequently, a meaningful validation requires a systematic, methodologically independent analysis of these approaches—a substantial research effort in its own right. Crucially, performing only limited testing in this context could yield misleading results; therefore, a dedicated study is necessary to properly isolate and quantify the importance of each stabilization term. However, the application of both terms is mandated despite their relative significance, as the weak form derived here requires them.
	
	A logical progression of this research would be a controlled inter‑model comparison across carefully selected parameter regimes. This study should evaluate solutions obtained via three methods at identical physical timestamps: an unstabilised approach using sufficiently small time steps, the stabilization method proposed by \citet{kausStabilizationAlgorithmGeodynamic2010}, and the full variational stabilization developed here. Such a benchmark is necessary to quantitatively assess where the Full stabilizationilization offers measurable advantages. Furthermore, to avoid misleading conclusions, these comparisons must span a range of settings rather than relying on a single benchmark.
	
	The approach highlights the fundamental difference between the free‑surface boundary and the traditionally used Neumann boundary condition for the free surface. We show that the presence of the free surface turns the problem into a time‑dependent equation even when inertial effects are ignored, in contrast to the steady‑state formulation commonly employed. We think that complementing the previous work on the stabilization term for free surfaces with the current mathematical foundation can pave the way for better time‑stepping methods that can improve accuracy and stability of models where surface and subsurface interaction is important.
	\appendix
	\pagebreak
	
	\section{Supplementary Material}
	\label{app:supplementary}
	
	This supplementary material provides additional details on the numerical
	workflow, Arbitrary Lagrangian--Eulerian (ALE) mesh movement, remeshing
	procedure, and verification benchmarks used to assess the numerical
	implementation. The main text introduces the governing equations and the
	physical model; the present section focuses on their numerical
	implementation and quantitative verification.
	
	\subsection{Numerical Workflow and ALE Time Integration}
	\label{subsec:ale_stokes_workflow}
	
	The numerical solution advances in time using an ALE formulation in which
	the physical velocity and the mesh velocity are treated as distinct
	quantities. At each time step, the Stokes problem is first solved on the
	current computational mesh. The resulting physical velocity is then used
	to determine the motion of selected mesh boundaries through a separate
	mesh-movement problem. The mesh coordinates are subsequently updated
	explicitly, and the resulting mesh quality is assessed before proceeding
	to the next time step.
	
	Let $t^n$ denote the current time and let $\Delta t^n$ be the time-step
	size. The computational state at time $t^n$ consists of the mesh
	coordinates, the physical fields, and the associated solver and function
	space definitions. The complete time-stepping procedure is summarized
	schematically in Algorithm~\ref{alg:ale-stokes}.
	
	\begin{algorithm}[H]
		\caption{ALE time-stepping procedure}
		\label{alg:ale-stokes}
		\begin{algorithmic}[1]
			\State Initialize mesh, physical fields, finite-element spaces,
			solvers, tolerances, and output schedule
			\For{$n=0,1,2,\ldots$ until $t^n\geq T$}
			\If{mesh or function-space reconstruction is required}
			\State Construct finite-element spaces
			$Z=U\times P$ and mesh space $U_m$
			\State Reinitialize the Stokes and mesh-movement solver
			\EndIf
			\State Assemble and solve the Stokes problem
			for $(\mathbf{u}^{n+1},p^{n+1})$
			\State Save solution and checkpoint data if scheduled
			\State Solve the mesh-movement problem for
			$\mathbf{u}_m^{n+1}$
			\State Update mesh coordinates according to
			$\mathbf{x}^{n+1}
			=
			\mathbf{x}^{n}
			+
			\Delta t^n\mathbf{u}_m^{n+1}$
			\State Evaluate the minimum mesh-quality measure
			$q_{\min}$
			\If{$q_{\min}<q_{\mathrm{thr}}$}
			\State Regenerate the computational mesh
			\EndIf
			\State Advance the time controller
			\EndFor
		\end{algorithmic}
	\end{algorithm}
	
	The separation between the physical and mesh velocities is important in
	the ALE formulation. The physical velocity $\mathbf{u}$ determines the
	evolution of the modeled material, whereas $\mathbf{u}_m$ controls the
	motion of the computational grid. This distinction allows the mesh to
	follow the evolving geometry while retaining sufficient control over
	element quality.
	
	\subsection{Mesh Movement}
	\label{subsec:mesh-movement}
	
	The mesh is moved by solving a separate linear elliptic smoothing problem
	for the mesh velocity $\mathbf{u}_m$. In the present implementation, the
	mesh velocity is obtained from
	\begin{equation}
		\int_{\Omega}
		\nabla\mathbf{u}_m:\nabla\mathbf{v}_m
		\,\mathrm{d}\Omega
		=
		0
		\qquad
		\forall\,\mathbf{v}_m\in U_m,
		\label{eq:mesh-movement}
	\end{equation}
	subject to Dirichlet conditions that prescribe the motion of boundaries
	that follow the material and constrain boundaries that must remain fixed.
	
	This harmonic-extension formulation propagates prescribed boundary
	motion into the interior of the domain while providing a smooth
	displacement field. The mesh coordinates are then advanced explicitly
	according to
	\begin{equation}
		\mathbf{x}^{n+1}
		=
		\mathbf{x}^{n}
		+
		\Delta t^n
		\mathbf{u}_m^{n+1}.
		\label{eq:ale-coordinate-update}
	\end{equation}
	
	The mesh-movement procedure is summarized in
	Algorithm~\ref{alg:mesh-movement}.
	
	\begin{algorithm}[H]
		\caption{Mesh-movement step}
		\label{alg:mesh-movement}
		\begin{algorithmic}[1]
			\State Assemble the elliptic variational problem for
			$\mathbf{u}_m$
			\State Apply mesh-motion and fixed-boundary conditions
			\State Solve the resulting linear system for $\mathbf{u}_m$
			\State Update the mesh coordinates using
			$\mathbf{x}^{n+1}
			=
			\mathbf{x}^{n}
			+
			\Delta t^n\mathbf{u}_m^{n+1}$
		\end{algorithmic}
	\end{algorithm}
	
	The mesh-movement solver is intentionally separated from the Stokes
	solver. Consequently, modifications to mesh motion or boundary tracking
	do not require changes to the physical flow formulation.
	
	\subsection{Mesh Quality and Remeshing}
	\label{subsec:mesh-remeshing}
	
	Mesh quality is monitored after each ALE update. The mesh generator stores
	the relevant boundary and interface markers together with a prescribed
	quality threshold $q_{\mathrm{thr}}$. A geometric quality measure is
	computed for every element; in the present implementation, the minimum
	element angle is used as the primary quality indicator.
	
	The global mesh-quality criterion is
	\begin{equation}
		q_{\min}
		=
		\min_{K\in\mathcal{T}_h} q_K,
		\qquad
		q_{\min}>q_{\mathrm{thr}},
		\label{eq:mesh-quality}
	\end{equation}
	where $q_K$ denotes the quality measure of element $K$ and
	$\mathcal{T}_h$ is the current mesh.
	
	When
	\begin{equation}
		q_{\min}<q_{\mathrm{thr}},
	\end{equation}
	the accumulated mesh distortion is considered excessive and the
	computational mesh is regenerated. The regenerated mesh is constructed
	from the current boundary and interface representations so that the
	prescribed geometry and physical interfaces are retained.
	
	Following remeshing, the finite-element spaces and associated solution
	fields are reconstructed on the new mesh, and the numerical solution is
	transferred to the regenerated discretization. This procedure prevents
	progressive element distortion from degrading the accuracy or robustness
	of the simulation.
	
	\subsection{Numerical Verification and Benchmarking}
	\label{subsec:numerical-verification}
	
	The numerical implementation was assessed using two complementary
	analytical benchmarks. The first employs a smooth manufactured solution
	and provides a controlled verification of the discretization, forcing,
	boundary conditions, and field evaluation. The second is the SolVi
	viscous-inclusion benchmark, which introduces a strong viscosity contrast
	and a discontinuous internal interface and therefore constitutes a more
	stringent test of the numerical formulation.
	
	For both benchmarks, the accuracy of the computed fields is quantified
	using
	\begin{align}
		E_{u,L^2}
		&=
		\left\|
		\mathbf{u}
		-
		\mathbf{u}^{\mathrm{exact}}
		\right\|_{L^2(\Omega)},
		\label{eq:velocity-l2-error}
		\\
		E_{p,L^2}
		&=
		\left\|
		p
		-
		p^{\mathrm{exact}}
		\right\|_{L^2(\Omega)}.
		\label{eq:pressure-l2-error}
	\end{align}
	
	Because the pressure in the incompressible Stokes problem is defined only
	up to an arbitrary additive constant, both the numerical and analytical
	pressures are normalized to zero mean prior to evaluating
	$E_{p,L^2}$. Specifically,
	\begin{equation}
		\bar{p}
		=
		\frac{1}{|\Omega|}
		\int_{\Omega}p\,\mathrm{d}\Omega,
		\qquad
		\widetilde{p}
		=
		p-\bar{p}.
		\label{eq:pressure-gauge}
	\end{equation}
	
	This normalization removes the arbitrary pressure gauge and ensures that
	the reported pressure error measures the spatial difference between the
	two fields.
	
	\subsection{Manufactured-Solution Benchmark}
	\label{subsec:manufactured-solution}
	This is a benchmark to demonstrate solver's ability to accurately solve a base Stokes flow problem. The benchmark is taken from \cite{buelerPETScPartialDifferential2021}. A smooth manufactured solution on the unit square,
	\[
	\Omega=(0,1)\times(0,1),
	\]
	is used as the baseline verification problem. The exact velocity and
	pressure fields are
	\begin{align}
		u_x^{\star}(x,y)
		&=
		\sin(4\pi x)\cos(4\pi y),
		\\
		u_y^{\star}(x,y)
		&=
		-\cos(4\pi x)\sin(4\pi y),
		\\
		p^{\star}(x,y)
		&=
		\pi\cos(4\pi x)\cos(4\pi y).
	\end{align}
	
	The forcing and boundary data are generated consistently from these
	analytical fields. The exact velocity is prescribed along the complete
	boundary:
	\begin{align}
		\mathbf{u}
		&=
		\left(0,-\sin(4\pi y)\right),
		&& x=0,1,
		\\
		\mathbf{u}
		&=
		\left(\sin(4\pi x),0\right),
		&& y=0,1.
	\end{align}
	
	The manufactured solution is smooth throughout the domain and contains
	no material discontinuity or internal interface. It therefore provides a
	controlled test of the implementation before more challenging
	heterogeneous configurations are considered.
	
	The computed errors for the reported discretization are
	\begin{table}[H]
		\centering
		\caption{
			Errors obtained for the manufactured-solution benchmark.
		}
		\label{tab:manufactured-errors}
		\begin{tabular}{lc}
			\hline
			Metric & Error \\
			\hline
			$E_{u,L^2}$
			& \texttt{2.330424e-04} \\
			$E_{p,L^2}$
			& \texttt{9.313305e-03} \\
			\hline
		\end{tabular}
	\end{table}
	
\begin{figure}[H]
		\centering
		\includegraphics[width=\textwidth]
		{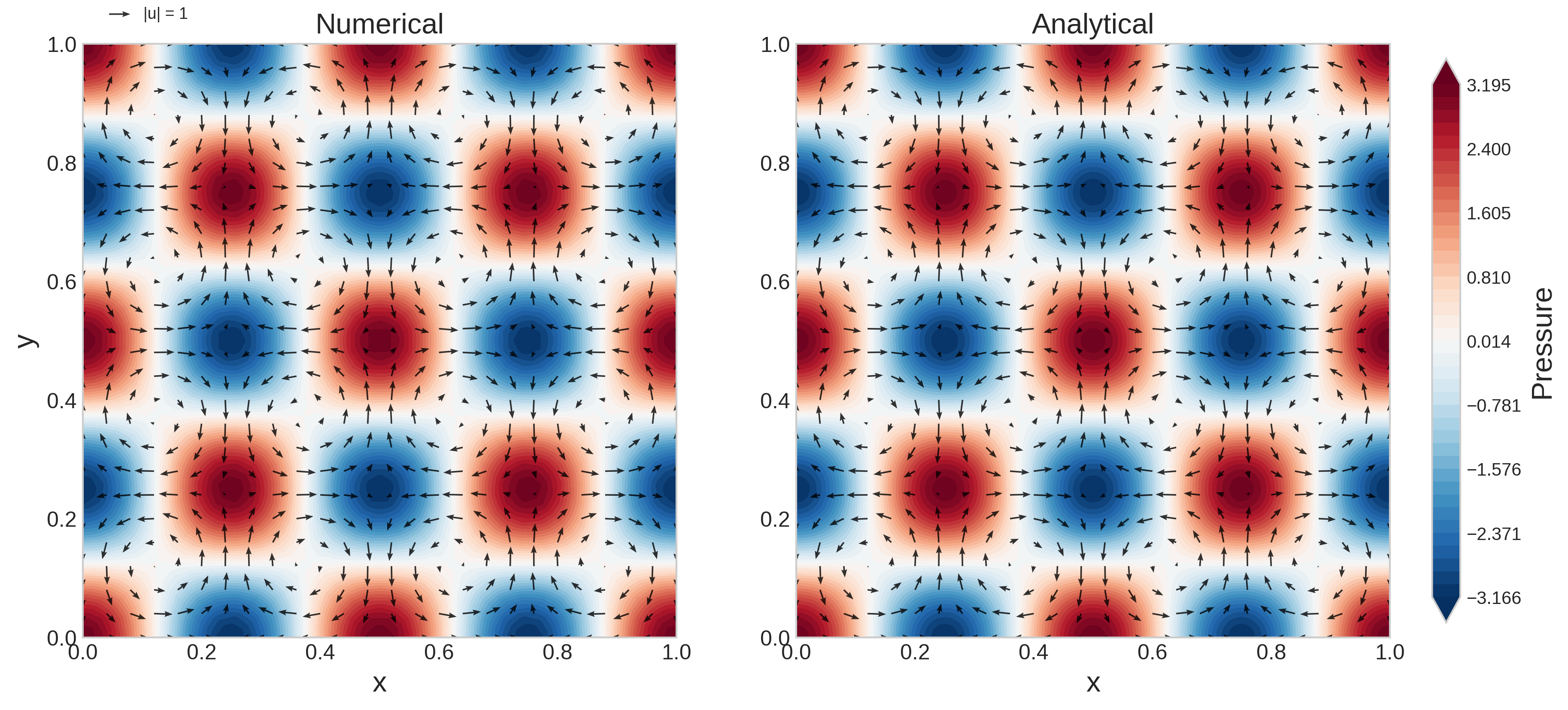}
		\caption{
			Manufactured-solution benchmark: numerical and analytical pressure
			fields with superimposed velocity vectors.
		}
		\label{fig:cavity-domain-compare}
	\end{figure}

\begin{figure}[H]
		\centering
		\includegraphics[width=0.78\textwidth]
		{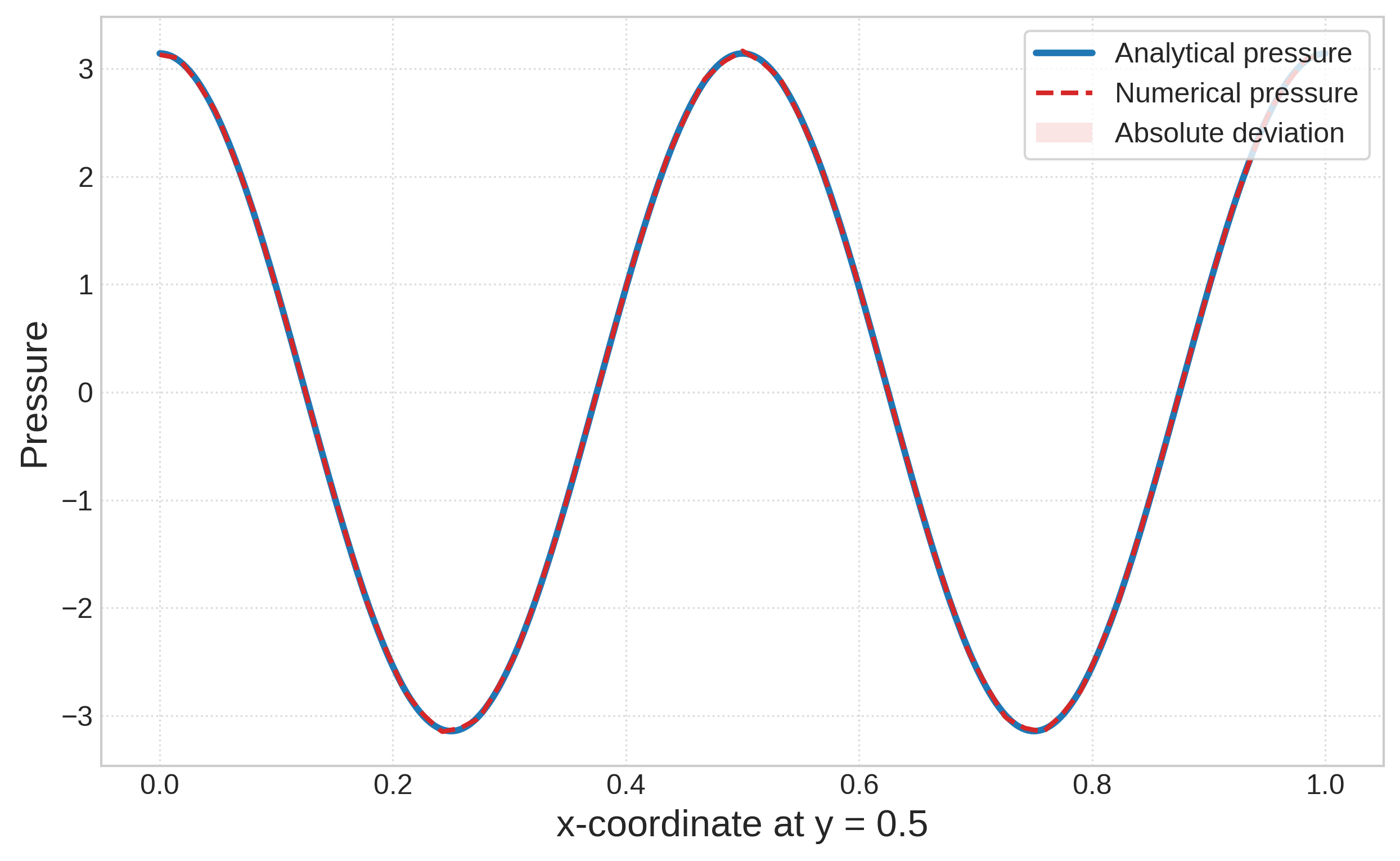}
		\caption{
			Manufactured-solution benchmark: numerical and analytical pressure
			along the horizontal centerline, $y=0.5$.
		}
		\label{fig:cavity-midline-pressure}
	\end{figure}

	The small velocity and pressure errors demonstrate that the numerical
	implementation reproduces the prescribed smooth solution accurately.
	Figure~\ref{fig:cavity-domain-compare} shows the spatial distributions
	of the numerical and analytical fields, while
	Figure~\ref{fig:cavity-midline-pressure} compares the corresponding
	pressure profiles along the horizontal centerline.
	
	\subsection{SolVi Viscous-Inclusion Benchmark}
	\label{subsec:solvi-benchmark}
	
	The second verification problem is the SolVi viscous-inclusion benchmark \cite{schmidAnalyticalSolutionsDeformable2003} that has been used extensively in the geodynamics community \cite{heisterHighAccuracyMantle2017,duretzFreeSurfaceCapturing2016}.
	In contrast to the smooth manufactured solution, this test contains a
	circular inclusion with a strong viscosity contrast relative to the
	surrounding matrix. It therefore evaluates the ability of the numerical
	implementation to represent a discontinuous material coefficient and the
	associated perturbation of the velocity and pressure fields.
	
	The inclusion is centered at the origin and has radius $R$. The viscosity
	field is defined by
	\begin{equation}
		\eta(x,y)
		=
		\begin{cases}
			\eta_i, & r<R,\\
			\eta_m, & r\geq R,
		\end{cases}
		\qquad
		r^2=x^2+y^2,
		\label{eq:solvi-viscosity}
	\end{equation}
	with the viscosity contrast parameter
	\begin{equation}
		B=
		\frac{\eta_i-\eta_m}
		{\eta_i+\eta_m}.
		\label{eq:solvi-B}
	\end{equation}
	
	For the sign convention adopted by the present implementation, the
	analytical velocity is
	\begin{align}
		u_x &=
		\begin{cases}
			-(1-B)\dot{\epsilon}\,x,
			& r<R,\\[4pt]
			-\dot{\epsilon}
			\left[
			x
			-\dfrac{BR^2x}{r^2}
			-\dfrac{BR^2(x^3-3xy^2)}{r^4}
			+\dfrac{BR^4(x^3-3xy^2)}{r^6}
			\right],
			& r\geq R,
		\end{cases}
		\\[6pt]
		u_y &=
		\begin{cases}
			+(1-B)\dot{\epsilon}\,y,
			& r<R,\\[4pt]
			-\dot{\epsilon}
			\left[
			-y
			+\dfrac{BR^2y}{r^2}
			-\dfrac{BR^2(3x^2y-y^3)}{r^4}
			+\dfrac{BR^4(3x^2y-y^3)}{r^6}
			\right],
			& r\geq R.
		\end{cases}
	\end{align}
	
	The corresponding pressure field is
	\begin{equation}
		p(x,y)
		=
		\begin{cases}
			0,
			& r<R,\\[4pt]
			-4\dot{\epsilon}\eta_m B R^2
			\dfrac{x^2-y^2}{r^4},
			& r\geq R.
		\end{cases}
		\label{eq:solvi-pressure}
	\end{equation}
	
	The analytical velocity is imposed on the external boundary, while the
	material contrast is represented directly through the discontinuous
	viscosity distribution in Eq.~\eqref{eq:solvi-viscosity}. Analytical and
	numerical pressures are normalized to zero mean before comparison.
	
	The resulting errors are summarized in
	Table~\ref{tab:solvi-error-norms}.
	
	\begin{table}[H]
		\centering
		\caption{
			Errors obtained for the SolVi viscous-inclusion benchmark after
			zero-mean pressure normalization.
		}
		\label{tab:solvi-error-norms}
		\begin{tabular}{lc}
			\hline
			Metric & Error \\
			\hline
			$E_{u,L^2}$
			& \texttt{7.005631e-04} \\
			$E_{p,L^2}$
			& \texttt{2.483134e-02} \\
			\hline
		\end{tabular}
	\end{table}
	
\begin{figure}[H]
		\centering
		\includegraphics[width=0.78\textwidth]
		{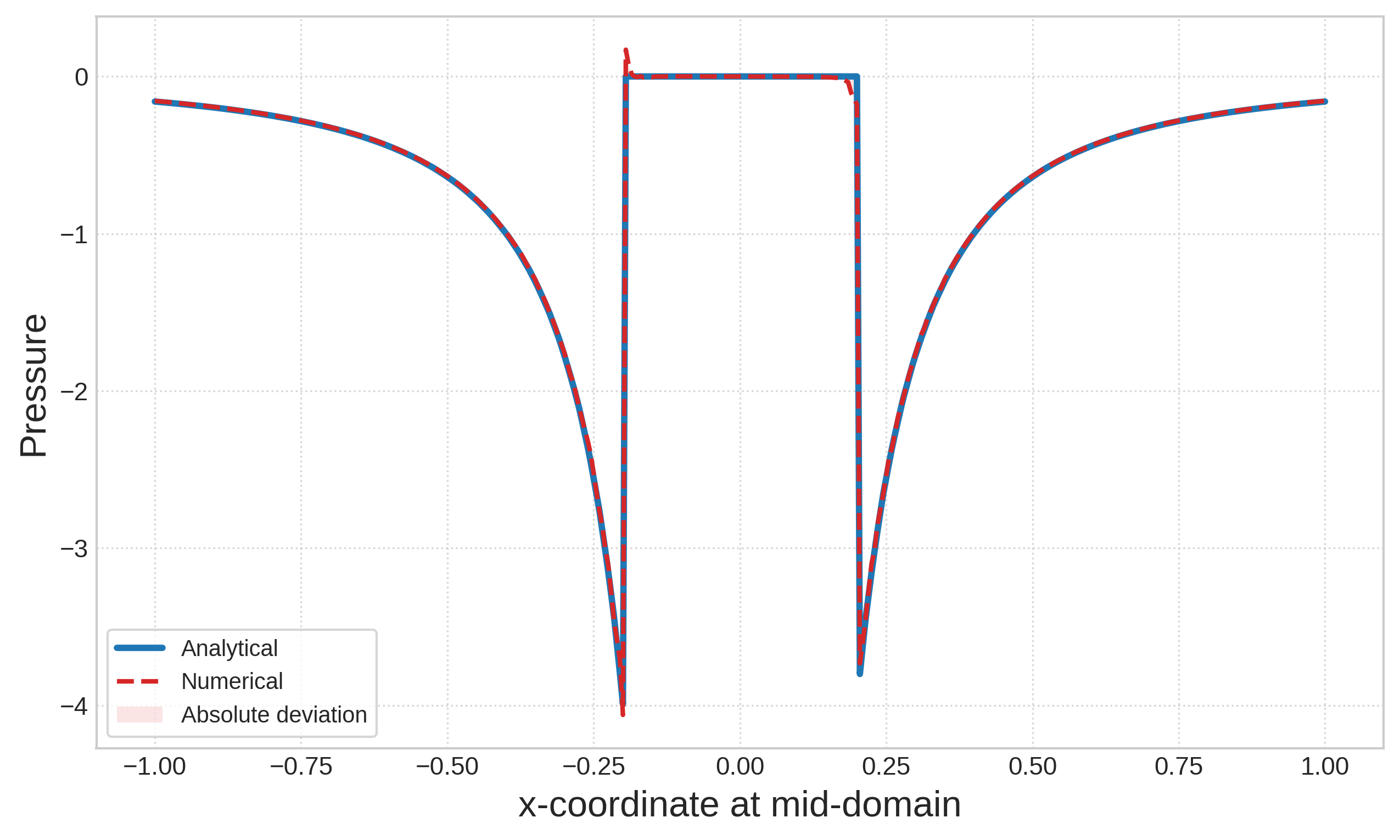}
		\caption{
			SolVi benchmark: numerical and analytical pressure along the
			horizontal centerline.
		}
		\label{fig:solvi-midline-pressure}
	\end{figure}

\begin{figure}[H]
		\centering
		\includegraphics[width=\textwidth]
		{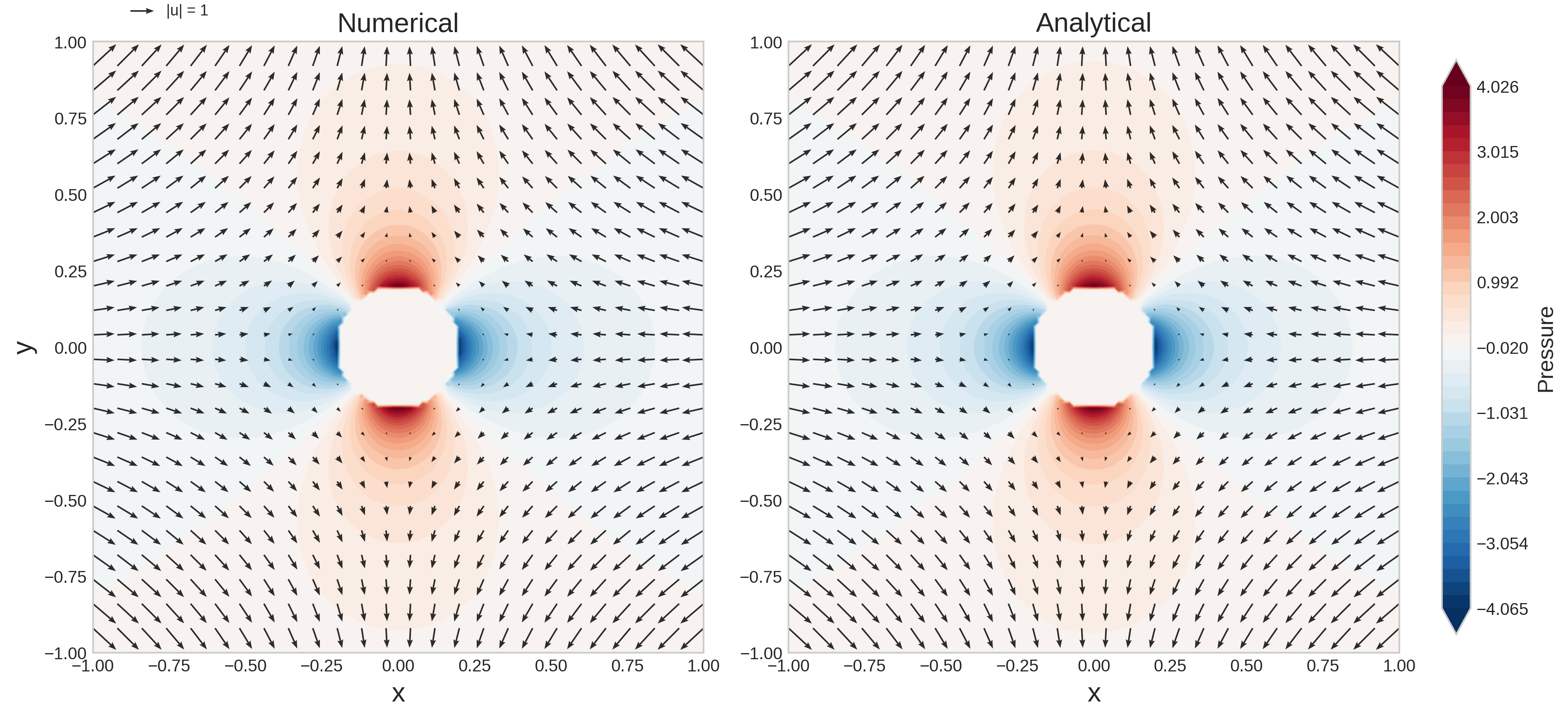}
		\caption{
			SolVi benchmark: comparison of numerical and analytical pressure
			fields with superimposed velocity vectors.
		}
		\label{fig:solvi-domain-overlay}
	\end{figure}

	Figure~\ref{fig:solvi-midline-pressure} compares the analytical and
	computed pressure along the horizontal centerline. The principal
	pressure extrema and the characteristic symmetry of the analytical
	solution are reproduced by the numerical solution. The corresponding
	domain-wide comparison is shown in
	Figure~\ref{fig:solvi-domain-overlay}, where the pressure field and
	velocity vectors illustrate the perturbation generated by the
	high-viscosity inclusion.
	
	The largest discrepancies occur in the vicinity of the circular
	viscosity interface. This behavior is expected because the material
	coefficient is discontinuous at $r=R$, producing a reduced regularity of
	the solution relative to the smooth manufactured benchmark.
	
	\subsection{Convergence and Verification Procedure}
	\label{subsec:convergence-verification}
	
	The two analytical benchmarks serve different but complementary purposes.
	The manufactured-solution benchmark provides a controlled test in which
	the analytical solution is smooth and known throughout the domain. It
	therefore establishes the baseline accuracy of the discretization and
	implementation. The SolVi benchmark adds a strong viscosity contrast and
	an internal discontinuous interface, providing a more demanding test of
	the formulation and its treatment of heterogeneous materials.
	
	The two tests therefore constitute a progressive verification framework:
	the manufactured solution establishes baseline numerical accuracy, while
	the SolVi benchmark evaluates the robustness of the implementation in
	the presence of strong material heterogeneity 
	and a sharp and well defined  internal interface. Together with the ALE mesh-movement and remeshing procedures
	described above, these benchmarks provide a reproducible framework for
	assessing the accuracy and numerical reliability of the simulation
	software.

	\bibliographystyle{elsarticle-num-names}
	
	\bibliography{references}
\end{document}